\documentclass[final]{siamltex}
\usepackage{xcolor}
\usepackage{amsmath}
\usepackage{amssymb}
\usepackage{mathrsfs}
\usepackage{graphicx}
\usepackage{bm}
\usepackage{ucs}
\usepackage[utf8x]{inputenc}
\usepackage{wrapfig}
\usepackage{color}
\usepackage{float}
\usepackage{epstopdf}
\usepackage{algorithm,algorithmic}
\usepackage[breaklinks,colorlinks=true,linkcolor=blue,citecolor=red,
  backref=page]{hyperref}

\DeclareMathAlphabet{\itbf}{OML}{cmm}{b}{it}
 \DeclareMathAlphabet\mathbfcal{OMS}{cmsy}{b}{n}

\renewcommand{\hat}{\widehat}
\def\EE{\mathbb{E}}
\def\RR{\mathbb{R}}

\def\bx{{{\itbf x}}}

\def\bM{{\itbf M}}
\def\bH{{\itbf H}}
\newcommand{\la}{\lambda}

\newcommand{\ep}{\varepsilon}

\newcommand{\om}{\omega}
\newcommand{\btau}{{\bm{\tau}}}
\newcommand{\bn}{{{\itbf n}}}
\newcommand{\cR}{\mathcal{R}}

\newtheorem{remark}[theorem]{Remark}

\begin{document}

\title{Pulse reflection in a random waveguide with a
  turning point} \author{Liliana Borcea\footnotemark[1] \and Josselin
  Garnier\footnotemark[2]} 
  
\maketitle


\renewcommand{\thefootnote}{\fnsymbol{footnote}}
\footnotetext[1]{Department of Mathematics, University of Michigan,
  Ann Arbor, MI 48109. {\tt borcea@umich.edu}}
\footnotetext[2]{Centre de Math\'ematiques Appliqu\'ees, Ecole Polytechnique, 91128 Palaiseau Cedex, France.  {\tt
    josselin.garnier@polytechnique.edu}}
\markboth{L. BORCEA, J. GARNIER}{Turning waves}

\begin{abstract}
We present an analysis of wave propagation and reflection in an acoustic waveguide
with random sound soft boundary and a turning point.  The waveguide
has slowly bending axis and variable cross section. The variation
consists of a slow and monotone change of the width of the waveguide
and small and rapid fluctuations of the boundary, on the scale of the
wavelength. These fluctuations are modeled as random. The turning
point is many wavelengths away from the source, which emits a pulse
that propagates toward the turning point, where it is reflected. To
focus attention on this reflection, we assume that the waveguide
supports a single propagating mode from the source to the turning
point, beyond which all the waves are evanescent.  We consider a
scaling regime where scattering at the random boundary has a
significant effect on the reflected pulse. In this regime 
scattering from the random boundary away from the turning point is negligible, 
while scattering from the random boundary around the turning point 
results in a strong, deterministic pulse deformation.
 The reflected pulse shape is not the same as the emitted one. It is
damped, due to scattering at the boundary, and is deformed by
dispersion in the waveguide. The reflected pulse also carries a random
phase.
\end{abstract}

\section{Introduction}
Guided waves arise in a wide range of applications in electromagnetics
\cite{collin1960field}, optics and communications
\cite{marcuse2013theory}, underwater acoustics \cite{kohler77}, and so
on. The classical theory of guided waves relies on the separability of
the wave equation in ideal waveguides with straight walls and filled
with homogeneous media \cite{lord1945theory}.  It decomposes the wave
field in independent waveguide modes, which are special solutions of
the wave equation. The modes are either propagating waves along the
axis of the waveguide or evanescent waves. They do not interact with
each other and have constant amplitudes determined by the source
excitation.

We study sound waves in two-dimensional waveguides with varying
cross section and slowly bending axis, where the waveguide effect is
due to reflecting boundaries, modeled for simplicity as sound
soft. The three-dimensional case and other boundary conditions can be
treated similarly, and do not involve conceptual differences.  We
refer to \cite{hazard2008improved,maurel2014propagation} for examples
of numerical studies of waves in slowly varying waveguides, and to
\cite{marcuse2013theory} for local mode decompositions of the wave
field, where the modes are coupled, and their amplitudes vary along
the waveguide axis.  An analysis of such a decomposition is
given in \cite{ahluwalia1974asymptotic,ting1983wave}, and the
transition of propagating modes to evanescent ones at turning points
in slowly changing waveguides is studied in
\cite{anyanwu1978asymptotic}.  Here we analyze this wave transition at
a turning point in a random waveguide with small and rapid random
fluctuations of the boundary on the scale of the wavelength, in
addition to the slow variations.

The wave field is generated by a source which emits a pulse with
central frequency $\om_o$ and bandwidth $B\ll \om_o$. It is the
superposition of a countable set of modes, of which only finitely many
propagate. To focus attention on the turning point, we consider a central frequency $\om_o$ 
such that there is a single propagating mode between the source and the turning point. We
also assume that the slow variation of the waveguide width is
monotone, so that no propagation occurs beyond the turning point.  Due
to energy conservation, the propagating mode is reflected at the
turning point and returns to the source location. The goal of the
paper is to analyze the pulse shape carried by this reflected wave.

Sound wave propagation in random waveguides is analyzed in
\cite{kohler77,dozier1978statistics,garnier_papa,garnier_single,gomez} for the case
of waveguides filled with a random medium, and in
\cite{alonso2011wave,borcea2014paraxial,gomez2011wave} for the case of
waveguides with random perturbations of straight boundaries.  We also
refer to \cite{marcuse2013theory,alonso2015electromagnetic} for the
analysis of electromagnetic waves in random waveguides.  A main
difficulty arising in the extension of these results to random
waveguides with slowly varying cross section is due to the turning
points, studied in this paper.

An analysis of random multiple scattering of turning waves is given in
\cite{kim1996stochastic,kim1996uniform}, in the context of wave
propagation in randomly layered media.  These results are relevant to
our study, specially the stochastic averaging theorem in
\cite{kim1996uniform}. In this paper we derive from first principles a
stochastic equation for the reflection coefficient of the propagating
mode in the random waveguide, and study in detail its statistics,
using the limit theorem in \cite{kim1996uniform}. To characterize the
reflected pulse, we carry out a multi-frequency analysis of the
reflection coefficient whose phase has a non-trivial random frequency dependence.
We quantify the standard deviation of the random fluctuations of the boundary 
{that trigger strong modifications of the amplitude and shape  of the reflected pulse.}
We show that such random 
fluctuations  have negligible effect on the pulse away from the turning point, but  near the turning point the effect is 
strong and leads to a deterministic pulse deformation and damping.
This pulse stabilization result is similar, but different from the ones
obtained in layered media in \cite{o1971reflections,clouet1994spreading,
  lewicki1994pulse,lewicki1994long}, in locally layered media in
\cite{solna2000ray}, in time-dependent layered media in
\cite{borcea2016pulse}, and in three-dimensional random media in
\cite{garnier2012coupled}. In these references the medium is random, not the boundary,
there is no turning point, and pulse deformation is observed when
 the standard deviation of the random fluctuations 
is larger than the one considered here. In this paper we explain why the random fluctuations have a stronger effect
close to the turning point than away from it.


The paper is organized as follows: We begin in section
\ref{sect:formulation} with the formulation of the problem and state
the pulse stabilization result in section \ref{sect:result}. The proof
of this result is in section \ref{sect:proof}.  We end with a summary
in section \ref{sect:sum}.
\section{Formulation of the problem}
\label{sect:formulation}
We describe in section \ref{sect:form1} the setup of the problem, and
define in section \ref{sect:form1p} the scaling regime. Then we give
in section \ref{sect:form2} the mode decomposition of the wave field,
and derive the stochastic differential equation satisfied by the
propagating mode. The remainder of the
paper is concerned with the analysis of this equation.

\subsection{Setup}
\label{sect:form1}

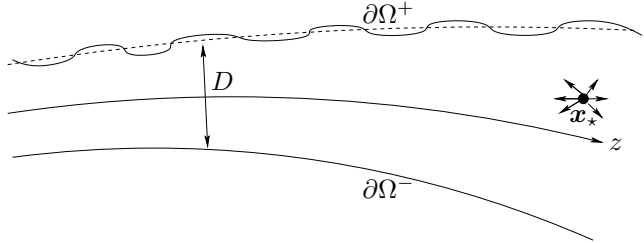
\begin{figure}[t]
\begin{picture}(0,0)%
\hspace{1in}\includegraphics[width = 0.65\textwidth]{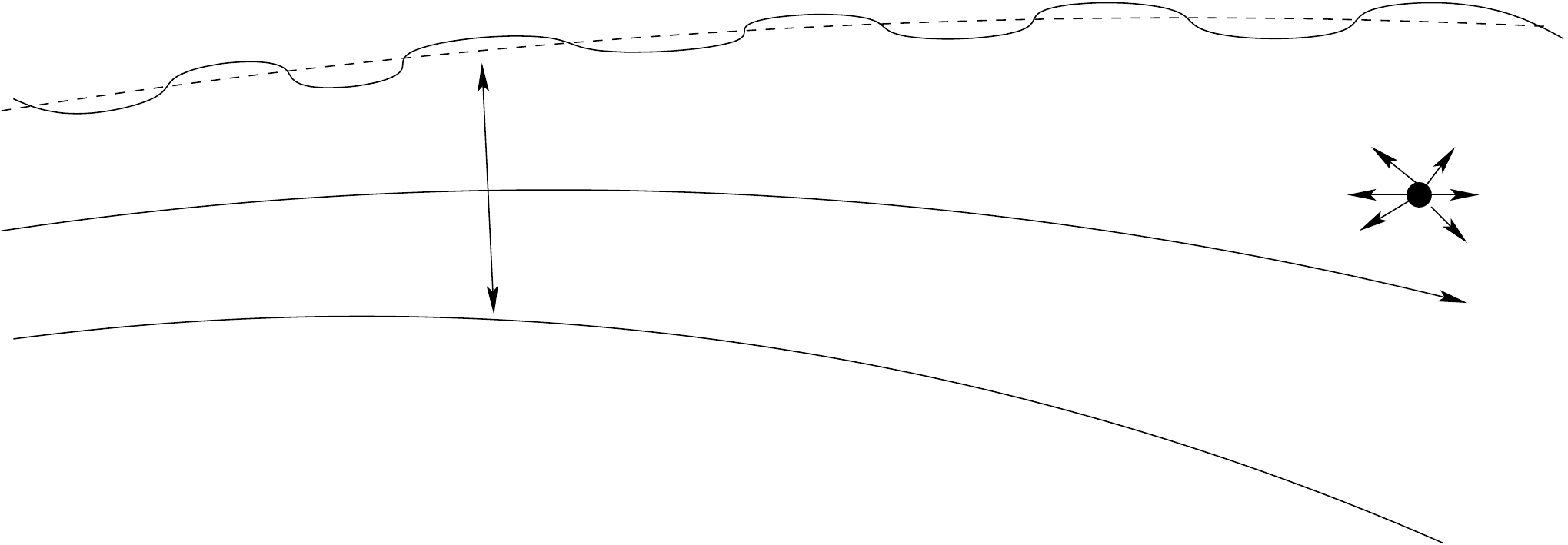}%
\end{picture}%
\setlength{\unitlength}{3947sp}%
\begingroup\makeatletter\ifx\SetFigFont\undefined%
\gdef\SetFigFont#1#2#3#4#5{%
  \reset@font\fontsize{#1}{#2pt}%
  \fontfamily{#3}\fontseries{#4}\fontshape{#5}%
  \selectfont}%
\fi\endgroup%
\begin{picture}(9770,1095)(1568,-1569)
\put(6550,-1000){\makebox(0,0)[lb]{\smash{{\SetFigFont{7}{8.4}{\familydefault}{\mddefault}{\updefault}{\color[rgb]{0,0,0}{\normalsize $z$}}%
}}}}
\put(6300,-820){\makebox(0,0)[lb]{\smash{{\SetFigFont{7}{8.4}{\familydefault}{\mddefault}{\updefault}{\color[rgb]{0,0,0}{\normalsize $\bx_\star$}}%
}}}}
\put(4050,-620){\makebox(0,0)[lb]{\smash{{\SetFigFont{7}{8.4}{\familydefault}{\mddefault}{\updefault}{\color[rgb]{0,0,0}{\normalsize $D$}}%
}}}}
\put(5000,-1300){\makebox(0,0)[lb]{\smash{{\SetFigFont{7}{8.4}{\familydefault}{\mddefault}{\updefault}{\color[rgb]{0,0,0}{\normalsize $\partial \Omega^-$}}%
}}}}
\put(5000,-200){\makebox(0,0)[lb]{\smash{{\SetFigFont{7}{8.4}{\familydefault}{\mddefault}{\updefault}{\color[rgb]{0,0,0}{\normalsize $\partial \Omega^+$}}%
}}}}
\end{picture}%
\caption{Illustration of a waveguide with monotonically increasing
  width $D$ and bending axis parametrized by the arc length $z$.
  The boundary $\partial \Omega$ is the union of the curves $\partial
  \Omega^-$ (the bottom boundary) and $\partial \Omega^+$ (the top
  boundary).  The top boundary is perturbed by small fluctuations
  modeled with a random process. The source of waves is at
  $\bx_\star$.  The waves first propagate towards negative $z$ in the form of a left-going propagating mode, 
  they are reflected at the turning point, and they propagate back towards positive $z$ in the form of a right-going mode. }
\label{fig:setup}
\end{figure}
Consider a two-dimensional waveguide occupying the semi-infinite
domain $\Omega$, with sound soft boundary $\partial \Omega = \partial
\Omega^- \cup \partial \Omega^+$ consisting of the union of two
curves, as illustrated in Figure \ref{fig:setup}. We refer to
$\partial \Omega^-$ as the bottom boundary and to $\partial \Omega^+$
as the top boundary.  The waveguide has a slowly bending axis
parametrized by the arc length $z$. Ideally, $\partial \Omega^-$ and
$\partial \Omega^+$ would be symmetric with respect to this axis, but
the top boundary is perturbed by small fluctuations.  The waveguide is
filled with a homogeneous medium with wave speed $c$, and the
excitation is due to a point source at location $\bx_\star \in
\Omega$, that emits the pulse
\begin{equation}
f(t) = \cos(\om_o t) F(Bt).
\label{eq:fm1}
\end{equation}
This pulse is modeled by a periodic carrier signal at frequency
$\om_o$, and a real-valued, smooth envelope function $F$ of dimensionless
argument.  Its Fourier transform $\hat F$ is
supported in the interval $[-\pi,\pi]$, so the Fourier transform
of \eqref{eq:fm1},
\begin{equation}
\hat f(\om) = \int_{-\infty}^\infty d t\, \cos(\om_o t) F(Bt) e^{i
  \om t} = \frac{1}{2B} \left[ \hat F\Big(\frac{\om-\om_o}{B}\Big) +
  \hat F\Big(\frac{\om+\om_o}{B}\Big) \right],
\label{eq:fm2}
\end{equation}
is supported in the frequency interval $[\om_o-\pi B, \om_o + \pi B]$
centered at $\om_o$, with bandwidth $B$, and its negative image
$[-\om_o-\pi B, -\om_o + \pi B]$. Since $F$ is real valued,
\begin{equation}
\hat F\Big(\frac{\om+\om_o}{B}\Big) = \overline{\hat
  F\Big(\frac{-\om-\om_o}{B}\Big)},
\label{eq:realval}
\end{equation}
where the bar denotes complex conjugate. We take $B \ll \om_o$, so
that $\la_o = {2 \pi c}/{\om_o}$ approximates the wavelength at all
frequencies $\om$ in the support of $\hat f(\om)$, and suppose that
$\la_o$ is small with respect to the arc length distance of order $L$
from the source to the turning point.

The wave field is modeled by the acoustic pressure $p(t,\bx)$, the
solution of the wave equation
\begin{equation}
\label{eq:waveeq}
\left(\Delta - \frac{1}{c^2} \partial_t^2\right) p(t,\bx) = f(t)
\delta(\bx-\bx_\star), \qquad \bx \in \Omega, ~ ~ t \in \RR, 
\end{equation}
with homogeneous Dirichlet boundary conditions 
\begin{equation}
p(t,\bx) = 0, \qquad \bx \in \partial \Omega, ~ ~ t \in \RR.
\label{eq:waveeq2}
\end{equation}
Prior to the excitation the medium is quiescent,
\begin{equation}
p(t,\bx) \equiv 0, \qquad t \ll 0.
\label{eq:waveeq1}
\end{equation}

It is convenient to write equations
\eqref{eq:waveeq}--\eqref{eq:waveeq2} in the orthogonal curvilinear
coordinate system with axes along $\btau(z/L)$ and $\bn(z/L)$, the
unit tangent and normal vectors to the axis of the waveguide, at arc
length $z$. These vectors change slowly in $z$, on the length scale $L
\gg \la_o$, according to the Frenet-Serret formulas
\begin{align}
\partial_z \btau\Big(\frac{z}{L}\Big) = \frac{1}{L}
\kappa\Big(\frac{z}{L}\Big) \bn\Big(\frac{z}{L}\Big), \qquad
\partial_z \bn\Big(\frac{z}{L}\Big) = - \frac{1}{L} \kappa
\Big(\frac{z}{L}\Big) \btau\Big(\frac{z}{L}\Big), \label{eq:fm3}
\end{align}
where $\kappa(z/L)$ is the curvature. We parametrize the 
points $\bx \in \Omega$ by $(r,z)$, using
\begin{equation}
\bx = \bx_{\parallel}(z) + r \bn\Big(\frac{z}{L}\Big),
\label{eq:fm4}
\end{equation}
where $\bx_{\parallel}(z)$ is on the waveguide axis, at arc length
$z$, and $r$ is the coordinate in the direction of the normal at
$z$. This coordinate lies in the interval $[r^-(z),r^+(z)]$, with
$r^-(z)$ at the bottom boundary $\partial \Omega^-$
\begin{equation}
r^-(z) = -\frac{D(z/L)}{2},
\label{eq:fm6p}
\end{equation}
and $r^+(z)$ at the randomly perturbed top boundary $\partial \Omega^+$
\begin{equation}
r^+(z)= \frac{D(z/L)}{2} \left[1 + 1_{(-\infty,0)}(z) \sigma \nu
  \Big(\frac{z}{\ell}\Big) \right].
\label{eq:fm6}
\end{equation}
Here $D(z/L)$ is the width of the unperturbed waveguide, a smooth (at least 
three times continuously differentiable) and monotonically increasing function that
varies slowly in $z$, on the scale $L$. The top boundary has small and
rapid random fluctuations on the left of the source, and $1_{(-\infty,0)}(z) $
is the indicator function of the negative axis $z < 0$, smoothed near
the origin.  The fluctuations are modeled by the zero-mean stationary
process $\nu$ of dimensionless argument, with autocorrelation function
\begin{equation}
\cR(\zeta) = \EE\big[ \nu(\zeta) \nu (0)\big].
\label{eq:autocorrel}
\end{equation}
This process is mixing, with rapidly decaying mixing rate, as defined
for example in \cite[section 2]{papanicolaou1974asymptotic}, and it is
bounded, with bounded first two derivatives, almost surely. We
normalize $\nu$ so that
\begin{equation}
\cR(0) = 1, \qquad \int_{-\infty}^\infty d \zeta \, \cR(\zeta) = 1 \mbox{ [or  $O(1)$]} ,
\end{equation} 
and control the amplitude of the fluctuations in \eqref{eq:fm6} by
the standard deviation $\sigma$, and their spatial scale  by the correlation length 
$\ell$.
The Fourier transform of $\cR$
$$
\widehat{\cR}(k) = \int_{-\infty}^\infty d\zeta \, e^{i k \zeta} \cR(\zeta) = \int_{-\infty}^\infty d\zeta \cos ( k \zeta) \cR(\zeta) 
$$
is the power spectral density of the stationary process $\nu$. It is an even and  nonnegative function.

In the curvilinear coordinate system the source is located at
$(r_\star, z=0)$, and the wave equation \eqref{eq:waveeq} becomes
\begin{align}
\left[\partial_r^2 -
  \frac{\frac{1}{L}\kappa\big(\frac{z}{L}\big)\partial_r}{1 -
    \frac{r}{L} \kappa\big(\frac{z}{L}\big)} +
  \frac{\partial_z^2}{\Big[ 1 - \frac{r}{L}
      \kappa\big(\frac{z}{L}\big)\Big]^2} + \frac{\frac{r}{L^2}
    \kappa'\big(\frac{z}{L}\big)\partial_z}{\Big[ 1 - \frac{r}{L}
      \kappa\big(\frac{z}{L}\big)\Big]^3} - \frac{1}{c^2} \partial_t^2
  \right] p(t,r,z) \nonumber \\= \left|1 - \frac{r_\star}{L}
\kappa(0)\right|^{-1} f(t)\delta(z) \delta(r-r_\star),\label{eq:fm7}
\end{align}
for $t \in \RR$, $z \in \RR$ and $r \in (r^-(z),r^+(z))$, with boundary
conditions \eqref{eq:waveeq2} given by
\begin{equation}
p(t,r^-(z),z) = p(t,r^+(z),z) = 0, \quad \forall \, t  \in \RR , ~  z \in \RR.
\label{eq:fm8}
\end{equation}
Here $\kappa'$ denotes the derivative of the curvature and we used the
parametrization \eqref{eq:fm4} of the points in the waveguide, with
\begin{equation*}
\partial_r \bx = \bn\Big(\frac{z}{L}\Big), \qquad \partial_z \bx =
\left[ 1 - \frac{r}{L} \kappa\Big(\frac{z}{L}\Big)
  \right]\btau\Big(\frac{z}{L}\Big),
\end{equation*}
and Lam\'{e} coefficients
\begin{equation*}
h_r = |\partial_r \bx| = 1, \qquad h_z = |\partial_z \bx| = \Big|1 - 
\frac{r}{L} \kappa\Big(\frac{z}{L}\Big)\Big|,
\end{equation*}
to write the  Laplacian
\begin{align*}
\Delta &= \frac{1}{h_r h_z} \left[ \partial_r \Big( \frac{h_r
    h_z}{h_r^2} \partial_r \Big) + \partial_z \Big( \frac{h_r
    h_z}{h_z^2} \partial_z \Big)\right],
\end{align*}
and the Dirac delta at $\bx_\star$, 
\[
\delta(\bx-\bx_\star) = \frac{1}{h_r h_z} \delta(z) \delta(r-r_\star).
\]
The problem is to analyze the wave field $p(t,r,z=0)$ at time $t >
T_f$, where $T_f$ is the duration of the emitted pulse $f(t)$. This
models the reflected wave in the random section of the waveguide,
which contains the turning point.

\subsection{The scaling regime}
\label{sect:form1p}
We define here a scaling regime where the random boundary fluctuations
have a significant effect on the reflected wave. The regime is defined
by the standard deviation $\sigma$ of the random fluctuations, and the
relation between the important length scales in the problem: the
central wavelength $\la_o = 2\pi c/\om_o$, the correlation length $\ell$ of the
random fluctuations, the scale $L$ of the slow variations of the
waveguide, and the width $D$ of the cross section.

The length scales are ordered as 
\begin{equation}
L \gg D \sim \la_o \sim \ell,
\label{eq:sc0}
\end{equation}
where $\sim$ denotes ``of the same order as".
In this scaling regime the central wavelength is of the same order as the correlation
length of the medium,  and is much smaller than the typical propagation distance,
  so that the waves interact efficiently with the boundary fluctuations. We model
\eqref{eq:sc0} using the small, dimensionless parameter
\begin{equation}
\ep = \frac{\ell}{L} \ll 1,
\label{eq:sc1}
\end{equation}
and use asymptotic analysis in the limit $\ep \to 0$ to characterize
the reflected wave.

The relation between the waveguide width $D(z/L)$ and the central
wavelength $\la_o$ determines the number
\[
N(z) = \lfloor {2 D(z/L)}/{\la_o} \rfloor
\]
of propagating modes in the local mode decomposition of the wave
$p(t,r,z)$, at given $z$, where $\lfloor ~ \rfloor$ denotes
the integer part. To simplify the analysis we assume that the central
frequency $\om_o$ of the pulse is such that $N(z) = 1$ for $z \in
(z_T(\om_o), 0)$, where $z_T(\om_o)<0$ is the arc length at the turning
point, satisfying
\begin{equation}
\la_o = 2 D\left({z_T(\om_o)}/{L}\right).
\label{eq:tp}
\end{equation}
The turning point is assumed simple, meaning that $D'(z_T(\om_o)/L) >
0$, and by the monotonicity of $D(z/L)$ we have $N(z) = 0$ for $ z <
z_T(\om_o)$.  Consistent with the slow variations of the waveguide on
the scale $L$, we suppose that $|z_T(\om_o)|  \sim L$.

We know from the study \cite{alonso2011wave} of waveguides with
randomly perturbed straight boundaries that the interaction of the
waves with the boundary fluctuations gives an order one net scattering
effect over the distance $L$ scaled as in \eqref{eq:sc1}, when the
standard deviation of the fluctuations is of the order
$\sqrt{\ep}$. Thus, we take
\begin{equation}
\sigma = \sqrt{\ep} \sigma_\ep,
\label{eq:sc2}
\end{equation}
with $\sigma_\ep$ at most of order one with respect to $\ep$. It will be adjusted later, so that
the effect of the random fluctuations of the waveguide boundary on the reflected pulse
is of order one as $\ep \to 0$. 

The duration $T_f$ of $f(t)$ is inverse
proportional to the bandwidth $B$, and must be much smaller than the
travel time from the source to the turning point and back, otherwise
$f(t)$ would not be a pulse. This implies the scaling relation 
\begin{equation}
\label{eq:sc3}
\frac{1}{\om_o} \ll \frac{1}{B} \ll \frac{L}{c} \sim  
\frac{1}{\ep  \om_o} ,
\end{equation}
where we used \eqref{eq:sc0} and $B \ll \om_o$. We show in section
\ref{sect:proof} that the characterization of the probability
distribution of the reflected pulse involves the joint distribution of
the reflection coefficients at  frequencies spaced by  $O(B)$.  We choose 
\begin{equation}
\frac{B}{\om_o} \sim \sqrt{\ep} ,
\label{eq:sc4}
\end{equation}
so that (\ref{eq:sc3}) is satisfied and the phases of the frequency-dependent reflection coefficients
have statistically dependent and independent components. 
This gives the pulse stabilization result after Fourier synthesis.

\subsubsection{The scaled variables}
We scale the arc length $z$ by $L$, and the waveguide width and
cross-range coordinate $r$ by $\ell$, 
\begin{equation}
\tilde z = z/L, ~ ~ \tilde r = r/\ell, ~ ~ \tilde D(\tilde z) =
D(z/L)/\ell.
\label{eq:sc5}
\end{equation}
The scaled frequency is
\begin{equation}
\tilde \om = \om \frac{\ell}{c} .
\label{eq:sc6}
\end{equation}
The central frequency is $\tilde\om_o = \om_o \ell / c =2 \pi \ell /\lambda_o$.
By (\ref{eq:sc4}) the bandwidth is such that $B \ell /c \sim \sqrt{\ep}$ so we introduce the scaled bandwidth $\tilde{B}$
defined by 
\begin{equation}
\tilde{B} = \frac{B \ell}{c} \frac{1}{\sqrt{\ep}} .
\end{equation}
The scaled wavenumber $k(\om) = \om/c$ is
\begin{equation}
\tilde k(\tilde\om) = k(\om) \ell .
\end{equation} 
All scaled quantities are of order one in the scaling regime described just above.

\subsubsection{The scaled equation}
We assume henceforth that the variables are scaled, and simplify the
notation by dropping the tilde. We take the Fourier transform of
\eqref{eq:fm7} with respect to time, and denote by 
$\hat p$ the wave field in the
scaled variables. After multiplying the
resulting equation by $L^2(1-\ep r \kappa(z))^2$ we obtain
\begin{align}
&\hspace{-0.15in}\left[\partial_z^2 + \frac{(1-\ep r
      \kappa(z))^2}{\ep^2}(k^2(\om)+\partial_r^2) - \frac{\kappa(z)(1-\ep r
      \kappa(z))}{\ep} \partial_r+ \frac{\ep r \kappa'(z) }{(1-\ep r
      \kappa(z))} \partial_z\right] \hat p(\om,r,z) \nonumber
  \\ &\hspace{3in}= \frac{\hat f^\ep(\om)}{\ep} \delta(r-r_\star)
  \delta(z),
\label{eq:md1}
\end{align} 
with 
\begin{equation}
\hat f^\ep(\om) = \frac{(1-\ep r_\star \kappa(0))}{ 2 
  \sqrt{\ep}  B} \left[ \hat F\Big(\frac{\om-\om_o}{\sqrt{\ep} B}\Big) +
  \overline{\hat F\Big(\frac{-\om-\om_o}{\sqrt{\ep} B}\Big)} \right], \label{eq:md2}
\end{equation}
and homogeneous Dirichlet boundary conditions 
\begin{equation}
\hat p(\om,r^\pm(z),z) = 0, \qquad r^-(z) = -\frac{D(z)}{2}, \quad
r^+(z) = \frac{D(z)}{2} \left[ 1 + \sqrt{\ep} \sigma_\ep
  \nu\Big(\frac{z}{\ep}\Big)\right],
\label{eq:md3} 
\end{equation}
for all $\om$ in the support of $\hat f^\ep(\om)$ and $z \in
\mathbb{R}$.  The wave is outgoing at $z > 0$, because there are no
random fluctuations there, and decays exponentially (is evanescent) at
$z < z_T(\om_o)$.

\subsection{Mode decomposition}
\label{sect:form2}
To define the mode decomposition, we change coordinates to map the
random boundary fluctuations to the coefficients of the wave equation
\eqref{eq:md1}. This way we obtain a linear differential operator
$\mathcal{L}^\ep$ that has an asymptotic expansion in $\ep$, and acts
on functions that vanish at the unperturbed boundary $r = \pm D(z)/2$
for all $z$. The modes are defined using the spectral decomposition of
the leading part of $\mathcal{L}^\ep$, and they have random amplitudes
satisfying stochastic differential equations driven by the process
$\nu$, with excitation given by jump conditions at $z=0$, where the
source lies.

\subsubsection{The random change of coordinates}
We use the following change of coordinates that maps the random boundary fluctuations to the 
wave operator
\begin{equation}
\label{eq:perteq1}
r = \rho + \frac{(2 \rho + D(z))}{4} \sqrt{\ep}
  \sigma_\ep\nu\Big(\frac{z}{\ep}\Big), \qquad \forall \, z < 0,
\end{equation}
where $\rho$ is in the unperturbed domain $[-D(z)/2,D(z)/2]$.  There are
no random fluctuations for $z >0$, so $r = \rho$ there.  Substituting
in \eqref{eq:md1} and using the chain rule, we obtain after
straightforward calculations that
\begin{equation}
\hat p^\ep(\om,\rho,z) = \hat p \Big(\om,\rho + \frac{(2 \rho +
  D(z))}{4} \sqrt{\ep} \sigma_\ep\nu\Big(\frac{z}{\ep}\Big),z \Big)
\end{equation}
satisfies the equation
\begin{equation}
\label{eq:perteq2}
\mathcal{L}^\ep \hat p^\ep(\om,\rho,z) = \frac{\hat f^\ep(\om)}{\ep}
\delta(\rho - r_\star) \delta(z),
\quad \forall \, \rho \in \left(-\frac{D(z)}{2}, \frac{D(z)}{2}
\right), ~ ~ z \in \mathbb{R},
\end{equation}
and the boundary conditions \eqref{eq:md3} become  
\begin{equation} \hat p^\ep(\om,\pm D(z)/2,z)
  = 0.
\end{equation}
The operator $\mathcal{L}^\ep$ is given by
\begin{align}
\label{eq:perteq3}\mathcal{L}^\ep = 
\sum_{j=0}^2 \ep^{j/2-2}\mathcal{L}_j + \partial_z^2 - \frac{(2
  \rho +
  D(z))}{2} \Big[\frac{\sigma_\ep}{\sqrt{\ep}}\nu'\Big(\frac{z}{\ep}\Big)
  - \frac{\sigma_\ep^2}{2}
  \nu\Big(\frac{z}{\ep}\Big)\nu'\Big(\frac{z}{\ep}\Big)
 \nonumber \\ + O(\sqrt{\ep})\Big]
\partial^2_{\rho z} + O(\ep)\partial_z +O(\ep^{-1/2}) \partial_\rho^2 + O(\ep^{-1/2}) \partial_\rho,
\end{align}
where $\mathcal{L}_j$ are differential operators with respect to
$\rho$, with coefficients that depend on $z$. These operators 
depend on $\ep$ only through $\sigma_\ep$ and the argument of $\nu$, $\nu'$, and $\nu''$. The leading operator
$\mathcal{L}_0$ is
\begin{align}
\mathcal{L}_0 &= k^2(\om) + \partial_\rho^2,
\end{align}
its first perturbation depends linearly on the random process $\nu$, 
\begin{align}
\mathcal{L}_1 = - \sigma_\ep  \nu
  \Big(\frac{z}{\ep}\Big) \partial_\rho^2  - \sigma_\ep \nu''
  \Big(\frac{z}{\ep}\Big)\frac{(2\rho + D(z))}{4} \partial_\rho,
\end{align}
the second perturbation is quadratic in $\nu$, 
\begin{align}
\mathcal{L}_2 = \frac{3\sigma_\ep^2}{4} \nu^2\Big(\frac{z}{\ep}\Big)
\partial_\rho^2 +
\sigma_\ep^2 {\nu'}^2\Big(\frac{z}{\ep}\Big)  \left[  \frac{(2
    \rho + D(z))^2}{16} \partial_\rho^2 + \frac{(2 \rho + D(z))}{4}
  \partial_\rho \right] \nonumber \\ +
\sigma_\ep^2\nu\Big(\frac{z}{\ep}\Big)\nu''\Big(\frac{z}{\ep}\Big)
\frac{(2 \rho + D(z))}{8} \partial_\rho- \kappa(z)\left[ 2 \rho
  (k^2(\om) + \partial_\rho^2) + \partial_\rho \right] .
\end{align}

\subsubsection{The waveguide modes}
The self-adjoint operator $\mathcal{L}_0$, acting on functions that
vanish at $\rho = \pm D(z)/2$ for any fixed $z$, has the
eigenfunctions
\begin{equation}
y_j(\rho,z) = \left[\frac{2}{D(z)}\right]^{1/2} \sin \left[\frac{(2
    \rho + D(z))}{2} \mu_j(z) \right], \quad \mu_j(z) = \frac{\pi
  j}{D(z)}, ~~ j = 1, 2, \ldots,
\label{eq:eigf}
\end{equation}
and eigenvalues $k^2(\om)-\mu_j^2(z)$, for $j = 1, 2, \ldots$ The
eigenfunctions form an orthonormal $L^2$ basis in $[-D(z)/2,D(z)/2]$,
so we can decompose the wave field at any $z$ as 
\begin{equation}
\hat p^\ep(\om,\rho,z) = \sum_{j=1}^\infty \hat u_j^\ep(\om,z)
y_j(\rho,z),
\label{eq:decomp}
\end{equation}
where $\hat u_j^\ep(\om,z)$ are waves in one dimension, called the
waveguide modes.  Substituting \eqref{eq:decomp} in
\eqref{eq:perteq2}, using the orthogonality of the eigenfunctions and
the identities given in appendix \ref{ap:identities}, we obtain
\begin{align}
&\hspace{-0.1in}\Big[ \partial_z^2 +
    \frac{k^2(\om)-\mu_j^2(z)}{\ep^2}\Big]\hat u_j^\ep(\om,z)
  +\frac{\sigma_\ep}{\ep^{3/2}} \left[\mu_j^2(z) \nu
    \Big(\frac{z}{\ep}\Big) + \frac{1}{4} \nu''
    \Big(\frac{z}{\ep}\Big)\right] \hat u_j^\ep(\om,z) \nonumber
  \\&\hspace{0.02in} + \frac{\sigma_\ep}{2 \ep^{1/2}}
  \nu'\Big(\frac{z}{\ep}\Big) \partial_z \hat u_j^\ep(\om,z) -
  \frac{\sigma_\ep^2}{\ep} \left\{ \frac{3 \mu_j^2(z)}{4}
  \nu^2\Big(\frac{z}{\ep}\Big) + \left[\frac{1}{8} + \frac{(\pi
      j)^2}{12}\right] {\nu'}^2\Big(\frac{z}{\ep}\Big) 
  \right. \nonumber \\ &\hspace{0.02in}\left. + \frac{1}{8}
  \nu\Big(\frac{z}{\ep}\Big) \nu''\Big(\frac{z}{\ep}\Big) \right\}\hat
  u_j^\ep(\om,z) - \frac{\sigma_\ep^2}{4} \nu\Big(\frac{z}{\ep}\Big)
  \nu'\Big(\frac{z}{\ep}\Big) \partial_z \hat u_j^\ep(\om,z) =
  \mathcal{C}_j^\ep\left(\om,z,\{\hat u_q^\ep\}_{q \ne j}\right),
\label{eq:decomp1}
\end{align}
for $z<0$. Here we neglected the remainder of order $\ep^{1/2}$, and
denoted by $\mathcal{C}_j^\ep$ the coupling terms that depend on the
modes $\hat u_q^\ep$, for $q \ne j$. The curvature $\kappa(z)$ of the
axis of the waveguide appears only in these terms.  The equations for
$z >0$ are simpler, because there are no random fluctuations in the
right-hand side. They, are obtained from \eqref{eq:decomp1} by setting
to zero all the terms that depend on the process $\nu$.

The first term in the  wave equations
\eqref{eq:decomp1} shows that the $j-$th mode is a propagating wave
when $\mu_j^2(z) < k^2(\om)$, and it is evanescent when the opposite
inequality holds.  By our scaling assumptions we have a single
propagating mode for $z<0$, the one indexed by $j = 1$. This interacts
with the evanescent modes via the coupling term
$\mathcal{C}_1^\ep$. We refer to \cite[section 3.3]{alonso2011wave}
for the analysis of such an interaction. It shows that the evanescent
modes can be expressed in terms of $u_1^\ep$, so that we can close the
wave equation for this propagating mode. We do not give here this
calculation, because it is basically the same as in
\cite{alonso2011wave}. The result is that the contribution of the
evanescent modes consists of an additional term in the equation for
$\hat u_1^\ep$, that is similar to the quadratic one in the
fluctuations, written in the curly bracket in \eqref{eq:decomp1}. We
will see in section \ref{sect:proof} that this term is negligible when
$\sigma_\ep$ is scaled so that the reflected pulse retains a
deterministic shape. For the sake of brevity, we do not include 
the contribution of the evanescent modes which play no role in the end.

\subsubsection{The equation for the propagating mode}
We can simplify the equation for the propagating mode $\hat u_1^\ep$
using integrating factors, by redefining the unknown
\begin{align}
\hspace{-0.0in} \hat u^\ep(\om,z) &= \hat u_1^\ep(\om,z) \exp
\hspace{-0.03in}\left[\frac{\ep^{1/2}\sigma_\ep}{4} \nu\Big(\frac{z}{\ep}\Big) -
  \frac{\ep \sigma_\ep^2}{16} \nu^2\Big(\frac{z}{\ep}\Big)\right] 
= \hat u_1^\ep(\om,z) \big[ 1 + O(\ep^{1/2})\big].
\label{eq:decomp3}\end{align}
Substituting in equation \eqref{eq:decomp1} for $j = 1$, we obtain
\begin{align}
\hspace{-0.05in}\partial_z^2 \hat u^\ep(\om,z) +
\left[\frac{k^2(\om)-\mu^2(z)}{\ep^2} + \frac{\sigma_\ep
    \mu^2(z)}{\ep^{3/2}} \nu\Big(\frac{z}{\ep}\Big) +
  \frac{\sigma_\ep^2}{\ep}
  g^\ep  (\om, z )\right] \hat
u^\ep(\om,z) = 0
\label{eq:decomp4}
\end{align}
for $z < 0$, with the simplified notation
\begin{align}
\mu(z) = \mu_1(z) = \frac{\pi}{D(z)}, \qquad 
g^\ep  (\om,z ) = -\frac{3}{4} \mu^2(z) \nu^2 \left( \frac{z}{\ep} \right) -
 \frac{\pi^2}{12}{\nu'}^2\left(\frac{z}{\ep} \right) 
 ,
\end{align}
where the contribution of the evanescent waves is not written as it vanishes in our scaling regime.
The excitation comes from the jump conditions at the source, with  $\hat f^\ep(\om)$ defined in \eqref{eq:md2},
\begin{align}
\hat u^\ep(\om,0^+)-\hat u^\ep(\om,0^-) &= 0, \label{eq:jump1}
\\ \partial_z \hat{u}^\ep(\om,0^+)-\partial_z \hat{u}^\ep(\om,0^-) & = \ep^{-1}\hat
  f^\ep(\om)y_1(r_\star,0).
\label{eq:jump2}
\end{align}

The remainder of the paper is concerned with the analysis of the
solution of \eqref{eq:decomp4}, with initial condition defined by
\eqref{eq:jump1}--\eqref{eq:jump2}, outgoing condition at $z >0$ and
exponential decay beyond the turning point, where the mode is
evanescent.
\section{The reflection coefficient and statement of results}
\label{sect:result}
We begin in section \ref{sect:res1} with the decomposition of $\hat
u^\ep(\om,z)$ in forward and backward going waves. This allows us to
define the reflection coefficient in section \ref{sect:res2}, and then
state the pulse stabilization result in section \ref{sect:res3}. This
result is derived in section \ref{sect:proof} under the assumption
that the turning point $z_T(\om)$ of the mode $\hat u^\ep(\om,z)$ is
simple, for any $\om$ in the support of $\hat f^\ep(\om)$. The
frequency-dependent turning point $z_T(\om)$ is defined by  
\begin{equation}
k(\om) = \mu(z_T(\om)) = \frac{\pi}{D(z_T(\om))},
\label{eq:res1}
\end{equation}
and it is unique due to the monotonicity of $D(z)$.
\subsection{The forward and backward going waves}
\label{sect:res1}
Let us write equation \eqref{eq:decomp4} as a first-order system of
stochastic differential equations
\begin{align}
\partial_z \begin{pmatrix} \hat u^\ep(\om,z)\\\hat
  v^\ep(\om,z)\end{pmatrix} &= \frac{i}{\ep} \begin{pmatrix} 0 & 1
  \\ k^2(\om) - \mu^2(z) &0 \end{pmatrix} \begin{pmatrix} \hat
  u^\ep(\om,z)\\\hat v^\ep(\om,z)\end{pmatrix} \nonumber \\ &+
\left[\frac{i \sigma_\ep}{\sqrt{\ep}} \mu^2(z) \nu
  \Big(\frac{z}{\ep}\Big) + i \sigma_\ep^2 g^\ep(\om, z) \right] \begin{pmatrix}0 & 0 \\ 1 &
  0 \end{pmatrix}
\begin{pmatrix} \hat
  u^\ep(\om,z)\\\hat v^\ep(\om,z)\end{pmatrix},
\label{eq:res3}
\end{align}
for the vector with components $\hat u^\ep(\om,z)$ and $ \hat
v^\ep(\om,z) = - i \ep \partial_z \hat u^\ep(\om,z).  $ Let also 
$\bM^\ep(\om,z)$ be a flow of smooth and invertible matrices, and 
define the vector
\begin{equation}
\begin{pmatrix} \hat a^\ep(\om,z) \\ \hat b^\ep(\om,z) \end{pmatrix} = 
\bM^{\ep,-1}(\om,z) \begin{pmatrix} \hat u^\ep(\om,z)\\\hat
  v^\ep(\om,z)\end{pmatrix},
\label{eq:res4}
\end{equation}
which satisfies equations 
\begin{align}
 & \partial_z \begin{pmatrix} \hat a^\ep(\om,z) \\ \hat
    b^\ep(\om,z) \end{pmatrix} = \bM^{\ep,-1}(\om,z) \left\{
  \frac{i}{\ep} \begin{pmatrix} 0 & 1 \\ k^2(\om) - \mu^2(z) &
    0 \end{pmatrix}\bM^\ep(\om,z)- \partial_z \bM^\ep(\om,z) \right.
  \nonumber \\ & \hspace{0.3in} \left. +\left[\frac{i
      \sigma_\ep}{\sqrt{\ep}} \mu^2(z) \nu \Big(\frac{z}{\ep}\Big) + i
    \sigma_\ep^2 g^\ep  (\om,z)\right]
 \begin{pmatrix}0 & 0 \\ 1 &
  0 \end{pmatrix} \bM^\ep(\om,z)\right\}\begin{pmatrix} \hat a^\ep(\om,z)
   \\ \hat b^\ep(\om,z) \end{pmatrix},
\label{eq:res5}
\end{align}
derived from \eqref{eq:res3}, where $\bM^{\ep,-1}$ denotes the inverse
of $\bM^\ep$.  The purpose of the decomposition \eqref{eq:res4} is to
remove the leading deterministic coupling term in \eqref{eq:res5} by a
proper choice of $\bM^\ep(\om,z)$, so that we can analyze the effect
of the random fluctuations. Then, we can associate the random 
fields $\hat a^\ep(\om,z)$ and $\hat b^\ep(\om,z)$ to the amplitudes
of the forward and backward going waves for the mode $\hat
u^\ep(\om,z)$, at $z > z_T(\om)$.

\subsubsection{The propagator}
The leading coupling term in \eqref{eq:res5} vanishes when
$\bM^\ep(\om,z) = \bM_\star^\ep(\om,z)$, the exact propagator matrix
in the unperturbed, slowly changing waveguide. This is the solution of
the flow problem
\[
\partial_z \bM_\star^\ep(\om,z) = \frac{i}{\ep} \begin{pmatrix} 0 & 1
  \\ k^2(\om) - \mu^2(z) & 0 \end{pmatrix} \bM_\star^\ep(\om,z), \quad
z<0,
\]
with $\bM_\star^\ep(\om,z=0)$ chosen so that we have the usual wave
decomposition at $z=0$, as in a waveguide with straight boundaries.
We work with an approximate propagator, which does not make the first
line in the right-hand side of \eqref{eq:res5} exactly zero, but it ensures that its
contribution to \eqref{eq:res5} converges to zero in the limit $\ep
\to 0$, uniformly in $z$, and its expression is explicit.

As in \cite{lynn1970uniform}, $\bM^\ep(\om,z)$ is the WKB
approximation of $\bM^\ep_\star(\om,z)$. It is a matrix with structure 
\begin{equation}
\bM^\ep(\om,z) = \begin{pmatrix} M_{11}^\ep(\om,z) & -
  \overline{M_{11}^\ep(\om,z)} \\ M_{21}^\ep(\om,z) &
  \overline{M_{21}^\ep(\om,z)} \end{pmatrix}, \qquad z < 0,
\label{eq:res6}
\end{equation}
where we recall that the bar denotes complex conjugate. The structure
in \eqref{eq:res6}
is like in waveguides with straight boundaries, and ensures energy
conservation, as follows later in the section. The entries in
\eqref{eq:res6} are defined in terms of the function
\begin{equation}
\phi_\om(z) =\left\{ \begin{aligned} &\int_{z_T(\om)}^z dz' \,
  \sqrt{k^2(\om)-\mu^2(z')}, \quad \qquad  ~z_T(\om) \le z \le 0,
  \\ &-\int^{z_T(\om)}_z dz' \, \sqrt{\mu^2(z')-k^2(\om)}, \qquad  z <
  z_T(\om),
\end{aligned} \right.
\label{eq:res7}
\end{equation}
which in turn defines
\begin{equation}
\eta_\om^\ep(z) = \left\{ \begin{aligned} & {\ep^{-2/3}} \left[ 3
    \phi_\om(z)/2\right]^{2/3}, \qquad \qquad z_T(\om) \le z \le 0,
  \\ & -\ep^{-2/3} \left[- 3 \phi_\om(z)/2\right]^{2/3}, \qquad  z <
  z_T(\om),
\end{aligned} \right.
\label{eq:res8}
\end{equation}
and 
\begin{equation}
Q_\om(z) = \left\{ \begin{aligned} & \frac{\left[ 3
      \phi_\om(z)/2\right]^{1/6}}{\left[k^2(\om)-\mu^2(z)\right]^{1/4}},
  \qquad z_T(\om) \le z \le 0, \\ & \frac{\left[ -3
      \phi_\om(z)/2\right]^{1/6}}{\left[\mu^2(z)-k^2(\om)\right]^{1/4}},
  \qquad z <z_T(\om).
\end{aligned} \right.
\label{eq:res9}
\end{equation}
Note that $Q_\om(z)$ is positive, at least twice continuously
differentiable, and at the turning point it satisfies
\begin{align}
Q_\om(z_T(\om)) = \gamma_\om^{-1/6}, \label{eq:res10} \quad
\partial_zQ_\om(z_T(\om)) = \frac{\theta_\om}{5\gamma_\om^{7/6} }, \quad 
\partial_z^2 Q_\om(z_T(\om)) = 
\frac{3 \rho_\om}{7 \gamma_\om^{7/6}} +
\frac{9 \theta_\om^2}{35 \gamma_\om^{13/6}},
\end{align}
where 
\begin{align}
\gamma_\om & = - \partial_z [\mu^2(z)]\Big|_{z=z_T(\om)} \hspace{-0.1in}=  \frac{2
  k^3(\om)}{\pi} D'(z_T(\om)) > 0, 
\label{eq:res13}\\
\theta_\om &= \frac{1}{2} \partial_z^2 [\mu^2(z)]\Big|_{z=z_T(\om)}, ~~~~
\rho_\om = \frac{1}{6} \partial_z^3 [\mu^2(z)]\Big|_{z=z_T(\om)} .
\end{align}
The function $\eta_\om^\ep(z)$ vanishes at the turning point, and its derivative
is given by
\begin{equation}
\partial_z \eta_\om^\ep(z) = \ep^{-2/3} Q^{-2}_\om(z), \qquad
\forall\, z < 0.
\label{eq:res14}
\end{equation}
The entries of the propagator matrix \eqref{eq:res6} are defined by
\begin{equation}
M_{11}^\ep(\om,z) = \ep^{-1/6}\sqrt{\pi} Q_\om(z) e^{-i\phi_\om(0)/\ep + i \pi/4} 
\left[A_i(-\eta_\om^\ep(z)) - i B_i(-\eta_\om^\ep(z))\right],
\label{eq:res15}
\end{equation}
and
\begin{align}
M_{21}^\ep(\om,z) =& - i \ep \partial_z M_{11}^\ep(\om,z) \nonumber
\\=& - \frac{\ep^{1/6} \sqrt{\pi}}{Q_\om(z)} e^{-i\phi_\om(0)/\ep - i
  \pi/4} \left[ A_i'(-\eta_\om^\ep(z)) - i
  B_i'(-\eta_\om^\ep(z))\right] \nonumber \\ &+ \ep^{5/6} \sqrt{\pi}
Q'_\om(z) e^{-i\phi_\om(0)/\ep - i \pi/4} \left[ A_i(-\eta_\om^\ep(z))
  - i B_i(-\eta_\om^\ep(z))\right],
\label{eq:res16}
\end{align}
in terms of the Airy functions \cite[chapter
  10]{abramowitz1972handbook} denoted by $A_i$ and $B_i$. 

The next lemma, proved in appendix \ref{ap:propagator}, shows that
$\bM^\ep(\om,z)$ approximates the exact propagator, and that
it is an invertible matrix with constant  determinant.

\vspace{0.05in}
\begin{lemma}
\label{lem.1}
The matrix-valued process \eqref{eq:res6}, with entries defined by
equations \eqref{eq:res15}--\eqref{eq:res16}, satisfies
\begin{equation}
\partial_z \bM^\ep(\om,z) = \frac{i}{\ep} \begin{pmatrix} 0 & 1
  \\ k^2(\om) - \mu^2(z) & 0 \end{pmatrix} \bM^\ep(\om,z)  -
\frac{i \ep Q_\om''(z)}{Q_\om(z)} \begin{pmatrix} 0 & 0 \\ 1 &
  0 \end{pmatrix} \bM^\ep(\om,z), 
\label{eq:res17}
\end{equation}
and 
\begin{equation}
\det \bM^\ep(\om,z) = 2, 
\label{eq:res18}
\end{equation}
for all $ z < 0$.
\end{lemma}

\vspace{0.05in} The next lemma, proved in appendix
\ref{ap:propagator}, describes the propagator as $z$ approaches $0$,
where the source lies.
\vspace{0.05in}
\begin{lemma}
\label{lem.2}
When $z<0$ and $|z| \ll 1,$ the entries
\eqref{eq:res15}--\eqref{eq:res16} of $\bM^\ep(\om,z)$ have the
following asymptotic expansions in $\ep$,
\begin{align*}
M_{11}^\ep(\om,z) = [k^2(\om)-\mu^2(0)]^{-1/4} \left\{ \exp \left[ \frac{i}{\ep} 
\big( \phi_\om(z)-\phi_\om(0)\big)\right] + O(\ep) \right\},
\end{align*}
and 
\begin{align*}
M_{21}^\ep(\om,z) = [k^2(\om)-\mu^2(0)]^{1/4} \left\{ \exp \left[
  \frac{i}{\ep} \big( \phi_\om(z)-\phi_\om(0)\big)\right] + O(\ep)
\right\}.
\end{align*}
The leading terms in these expansions are the entries of the
propagator in waveguides with straight boundaries and width $D(0)$.
\end{lemma}

\vspace{0.05in} \noindent Using this result in \eqref{eq:res4} and
\eqref{eq:res6}, we obtain a 
wave decomposition like  in waveguides with straight walls \cite[chapter
  20]{fouque07}. The wave field is, for $|z|\ll 1$,
\begin{align}
\hat u^\ep(\om,z) \approx [k^2(\om)-\mu^2(0)]^{-1/4} \left[ \hat
  a^\ep(\om,z) \exp \left( \frac{i}{\ep}\int_0^z dz' \,
  \sqrt{k^2(\om)-\mu^2(z')}\right) \right. \nonumber \\ \left. - \hat
  b^\ep(\om,z) \exp \left(- \frac{i}{\ep}\int_0^z dz' \,
  \sqrt{k^2(\om)-\mu^2(z')}\right) \right], \label{eq:res21}
\end{align}
and its derivative is 
\begin{align}
 \partial_z \hat u^\ep(\om,z) \approx \frac{i}{\ep}
         [k^2(\om)-\mu^2(0)]^{1/4} \left[ \hat a^\ep(\om,z) \exp
           \left( \frac{i}{\ep}\int_0^z dz' \,
           \sqrt{k^2(\om)-\mu^2(z')}\right) \right. \nonumber
           \\ \left. + \hat b^\ep(\om,z) \exp \left(-
           \frac{i}{\ep}\int_0^z dz' \,
           \sqrt{k^2(\om)-\mu^2(z')}\right) \right], \label{eq:res22}
\end{align}
with relative error of order $\ep$. 

In the vicinity of the turning point, for $|z - z_T(\om)| =
O(\ep^{2/3})$, the Airy functions and their derivatives are bounded,
as are $Q_\om(z)$ and its derivatives, described in \eqref{eq:res10}.  We obtain that the
entries in the first row of $\bM^\ep(\om,z)$ are large, of order
$\ep^{-1/6}$, and the entries in the second row are small, of order
$\ep^{1/6}$.

Beyond the turning point, at $z_T(\om)-z \gg O(\ep^{2/3})$, the
entries of $\bM^\ep(\om,z)$ grow exponentially, as stated in the next
lemma, proved in appendix \ref{ap:propagator}. The mode $\hat
u^\ep(\om,z)$ is evanescent in this region of the waveguide, and must
be exponentially decaying away from $z_T(\om)$.  This is ensured by
carefully chosen boundary conditions of the mode amplitudes, as
explained in the next section.

\vspace{0.05in}
\begin{lemma} 
\label{lem.3}
When $z_T(\om) - z \gg O(\ep^{2/3}),$ the entries
of the approximate propagator $\bM^\ep(\om,z)$ have the following
asymptotic expansions in $\ep$,
\begin{equation*}
M_{11}^\ep(\om,z) \approx [\mu^2(z)-k^2(\om)]^{-1/4} \exp \left[
  \frac{1}{\ep} \int_z^{z_T(\om)} \hspace{-0.1in} d z' \,
  \sqrt{\mu^2(z')-k^2(\om)} - \frac{i \phi_\om(0)}{\ep} - \frac{i
    \pi}{4} \right],
\end{equation*}
and 
\begin{equation*}
M_{21}^\ep(\om,z) \approx [\mu^2(z)-k^2(\om)]^{1/4} \exp \left[
  \frac{1}{\ep} \int_z^{z_T(\om)} \hspace{-0.1in} d z' \,
  \sqrt{\mu^2(z')-k^2(\om)} - \frac{i \phi_\om(0)}{\ep} + \frac{i
    \pi}{4} \right],
\end{equation*}
with relative error of order $\ep$. 
\end{lemma}

\vspace{0.05in} \noindent The Airy function $A_i$ and its derivative $A_i'$
decay exponentially in this region and are negligible. The asymptotic expansions above are determined by
$B_i$ and $B_i'$.

\subsubsection{The mode amplitudes}
To derive the system of differential equations satisfied by the mode
amplitudes, we use Lemma \ref{lem.1}  and the inverse of the propagator 
\begin{equation}
\bM^{\ep,-1}(\om,z) = \frac{1}{2} \begin{pmatrix}
  \overline{M_{21}^\ep(\om,z)} & \overline{M_{11}^\ep(\om,z)} \\
-M_{21}^\ep(\om,z) & M_{11}^\ep(\om,z) \end{pmatrix},
\label{eq:MInv}
\end{equation}
in equation \eqref{eq:res5}.  We obtain that
\begin{equation}
\partial_z\begin{pmatrix} \hat a^\ep(\om,z) \\ \hat
b^\ep(\om,z) \end{pmatrix} = \bH^\ep(\om,z) \begin{pmatrix} \hat
  a^\ep(\om,z) \\ \hat b^\ep(\om,z) \end{pmatrix}, \quad z < 0,
\label{eq:res23}
\end{equation}
with matrix-valued random process 
\begin{equation}
\bH^\ep(\om,z) = \begin{pmatrix} H_{11}^\ep(\om,z) &
  \overline{H_{21}^\ep(\om,z)} \\ H_{21}^\ep(\om,z) &
  -H_{11}^\ep(\om,z) \end{pmatrix},
\label{eq:res24}
\end{equation}
and entries equal, up to negligible terms, to
\begin{align} 
H_{11}^\ep(\om,z) &=\frac{1}{2}\left[ \frac{i \sigma_\ep}{\sqrt{\ep}}
  \mu^2(z) \nu \Big(\frac{z}{\ep}\Big) + i \sigma_\ep^2 g^\ep(\om,z) \right]
\Big|M_{11}^\ep(\om,z)\Big|^2, \label{eq:res25}\\ H_{21}^\ep(\om,z)
&=\frac{1}{2}\left[ \frac{i \sigma_\ep}{\sqrt{\ep}} \mu^2(z) \nu
  \Big(\frac{z}{\ep}\Big) + i \sigma_\ep^2 g^\ep(\om,z)\right]
\Big(M_{11}^\ep(\om,z)\Big)^2. \label{eq:res26}
\end{align}

To specify the solution of \eqref{eq:res23}, we need boundary
conditions. At $ z= 0$ we obtain from the jump conditions
\eqref{eq:jump1}--\eqref{eq:jump2}, equations
\eqref{eq:res21}--\eqref{eq:res22} and the outgoing condition $\hat
b^\ep(\om,0^+) = 0$ that
\begin{equation}
\hat b^\ep(\om,0^-) =  \frac{i C_F}{\sqrt{\ep} B} \left[ \hat
  F\Big(\frac{\om - \om_o}{\sqrt{\ep}B}\Big) + \overline{\hat
    F\Big(\frac{-\om - \om_o}{\sqrt{\ep}B}\Big)} \right],
\label{eq:res27}
\end{equation}
with constant
\begin{equation}
C_F = \frac{y_1(r_\star,0)}{4 [k^2(\om_o)-\mu^2(0)]^{1/4}}.
\label{eq:res28}
\end{equation}
Here we neglected the $O(\ep^{1/2})$ residual in the definition
\eqref{eq:md2} of $\hat f^\ep(\om)$, and in the expansion of $k(\om)$
for $\om = \om_o + O(\ep^{1/2})$. The forward going wave amplitudes at
$z = 0$ satisfy the relation
\begin{equation}
a^\ep(\om,0^-)-a^\ep(\om,0^+) = b^\ep(\om,0^-),
\label{eq:res29}
\end{equation}
and we need one more boundary condition.  This will ensure that the
wave is exponentially decaying away from the turning point.

The asymptotic expansion of the propagator at $z < z_T(\om)$, given in
Lemma \ref{lem.3}, shows that $\bM^\ep(\om,z)$ has exponentially
growing terms, due to the Airy function $B_i$. To compensate this
growth, we introduce here a boundary condition at some $z_{\rm b}$ far
enough from the turning point\footnote{We show in section
  \ref{sect:proof} that the result does not depend on the value
  of the fictitious boundary $z_{\rm b}$.}, satisfying $ z_T(\om) - z_{\rm b} \gg O(\ep^{2/3})$.
Definitions \eqref{eq:res4} and \eqref{eq:res15}--\eqref{eq:res16}
give that
\begin{align}
\hat u^\ep(\om,z_{\rm b}) = &\ep^{-1/6} C^\ep B_i(-\eta_\om^\ep(z_{\rm b}))
\left[a^\ep(\om,z_{\rm b}) - i e^{2 i \phi_\om(0)/\ep} b^\ep(\om,z_{\rm b})
  \right] \nonumber \\ &+ i \ep^{-1/6}C^\ep A_i(-\eta_\om^\ep(z_{\rm b}))
\left[ a^\ep(\om,z_{\rm b}) + i e^{2 i \phi_\om(0)/\ep}
  b^\ep(\om,z_{\rm b})\right],\label{eq:BC1}
\end{align}
and 
\begin{align}
\partial_z\hat u^\ep&(\om,z_{\rm b}) = -\ep^{-5/6} C^\ep\left[
  Q_{\om}^{-2}(z_{\rm b})B_i'(-\eta_\om^\ep(z_{\rm b})) - \ep
  \frac{Q_\om'(z_{\rm b})}{Q_\om(z_{\rm b})} B_i(-\eta_\om^\ep(z_{\rm b})) \right]
\nonumber \\ &\hspace{-0.05in}\times \left[a^\ep(\om,z_{\rm b}) - i e^{2 i
    \phi_\om(0)/\ep} b^\ep(\om,z_{\rm b}) \right] - i\ep^{-5/6} C^\ep
Q_{\om}^{-2}(z_{\rm b}) A_i'(-\eta_\om^\ep(z_{\rm b}))\nonumber
\\& \hspace{-0.05in}\times \left[ a^\ep(\om,z_{\rm b}) + i e^{2 i
    \phi_\om(0)/\ep} b^\ep(\om,z_{\rm b})\right] \left[1 +
  O(\ep)\right],\label{eq:BC2}
\end{align}
with constant
\[
C^\ep = \sqrt{\pi} Q_\om(z_{\rm b}) e^{-i\phi_\om(0)/\ep - i \pi/4}.
\]
We set to zero the coefficients of $B_i$ and $B_i'$ in these
expressions, to get an exponentially small wave field, and obtain the
boundary condition
\begin{equation}
a^\ep(\om,z_{\rm b}) = i e^{2 i \phi_\om(0)/\ep} b^\ep(\om,z_{\rm b}).
\label{eq:res33}
\end{equation}

\subsection{The reflection coefficient}
\label{sect:res2}
The mode amplitudes define the reflection coefficient
\begin{equation}
\hat R^\ep(\om,z) = \frac{\hat a^\ep(\om,z)}{\hat b^\ep(\om,z)},
\label{eq:res34}
\end{equation}
which is a complex number with modulus one. This is because 
the structure \eqref{eq:res24} of the matrix $\bH^\ep(\om,z)$ in
\eqref{eq:res23} ensures $\partial_z [ |\hat a^\ep(\om,z)|^2 - |\hat b^\ep(\om,z)|^2 ] = 0$, 
which gives the flux energy conservation equation
\begin{align}
|\hat a^\ep(\om,z)|^2 - |\hat b^\ep(\om,z)|^2 = \mbox{constant},
\qquad \forall \, z \in (z_{\rm b},0),
\end{align}
where the constant must equal zero by \eqref{eq:res33}. Thus,  we can write \eqref{eq:res34}  in the form
\begin{equation}
\hat R^\ep(\om,z) = i \exp \left[2 i \frac{\phi_\om(0)}{\ep} + i
  \psi_\om^\ep(z)\right],
\label{eq:res35} 
\end{equation}
with real-valued, random phase $\psi_\om^\ep(z)$. It satisfies the
differential equation
\begin{align}
\partial_z \psi_\om^\ep(z) = &\frac{2 \pi Q_\om^2(z)}{\ep^{1/3}} \Big[ \frac{
    \sigma_\ep}{\sqrt{\ep}} \mu^2(z) \nu \Big(\frac{z}{\ep}\Big) +
  \sigma_\ep^2 g^\ep(\om,z)\Big]
\big[A_i^2(-\eta_\om^\ep(z)) + B_i^2(-\eta_\om^\ep(z))\big] \nonumber
\\ &\times \cos^2 \Big\{ \frac{\psi_\om^\ep(z)}{2} -
\arg\left[A_i(-\eta_\om^\ep(z)) + i B_i(-\eta_\om^\ep(z))\right]\Big\},
\quad z > z_{\rm b},
\label{eq:res36ODE}
\end{align}
derived from \eqref{eq:res23} and \eqref{eq:res34}--\eqref{eq:res35},
with homogeneous boundary condition
\begin{equation}
\psi_\om^\ep(z_{\rm b}) = 0.
\label{eq:res36}
\end{equation}
We are particularly interested in the phase $\psi_\om^\ep(z = 0)$,
which defines the frequency-dependent amplitude of the reflected wave at the source.

\subsection{The pulse stabilization result}
\label{sect:res3}
Equations \eqref{eq:decomp}, \eqref{eq:res21}, \eqref{eq:res27},
\eqref{eq:res34} and \eqref{eq:res35} give
that the reflected pressure wave at $z = 0^-$ is given by
\begin{align}
\hat p_{\rm ref}^\ep(t,\rho,0^-) & = y_1(\rho,0)  \int \frac{d \om}{2 \pi} \,  \frac{1}{[k^2(\om) -
  \mu^2(0)]^{1/4}} e^{-i \om t/\ep} \hat a^\ep(\om,0^-)\nonumber \\
  &\approx -\frac{y_1(\rho,0)y_1(r_\star,0)}{2   \sqrt{k^2(\om_o)-\mu^2(0)}}f_{\rm ref}^\ep(t)\big[1+O(\sqrt{\ep})\big],
\label{eq:PS1}
\end{align}
with reflected pulse  
\begin{align*}
f_{\rm ref}^\ep(t) = \mbox{Re} \left\{\int \frac{d w}{2 \pi B} \hat F\Big(\frac{w}{B}\Big)
\exp\left[\frac{i [2\phi_{\om_o + \sqrt{\ep} w}(0) - (\om_o
      +\sqrt{\ep}w)t]}{\ep} + i \psi_{\om_o + \sqrt{\ep} w}^\ep(0)
  \right]\right\}.
\end{align*}
The wave emerging at the right of the source, at $z = 0^+$, 
is 
\begin{align}
\hat p_{\rm ref}^\ep(t,\rho,0^+) = y_1(\rho,0)  \int \frac{d \om}{2 \pi} \,\frac{1}{[k^2(\om) -
  \mu^2(0)]^{1/4}} e^{-i \om t/\ep} \hat a^\ep(\om,0^+),
\label{eq:PS1p}
\end{align}
with $a^\ep(\om,0^+)$ obtained from \eqref{eq:res29}. It is the
superposition of the reflected wave \eqref{eq:PS1} and the direct wave
that has no interaction with the random section of the waveguide, and
propagates from the source in the forward direction.

To describe $f_{\rm ref}^\ep(t)$ in the limit $\ep \to
0$, we eliminate first the large deterministic phase of the
integrand. For this purpose, we expand
\begin{equation*}
2\phi_{\om_o + \sqrt{\ep} w}(0) =2 \phi_{\om_o}(0) + \sqrt{\ep} w
T_{\om_o} + \ep w^2 \beta_{\om_o} + O(\ep^{3/2}),
\end{equation*}
where 
\begin{equation}
T_{\om_o} = 2 \partial_{\om} \phi_\om(0)\big|_{\om = \om_o} = \frac{2 k^2(\om_o)}{\om_o}
\int_{z_{T}(\om_o)}^0 \frac{d z }{\sqrt{k^2(\om_o)-\mu^2(z)}},
\end{equation}
is the travel time of the propagating mode from the source to the turning point
and back, and
\begin{align*}
\nonumber
\beta_{\om_o} =& \partial_\om^2 \phi_{\om}(0) \big|_{\om = \om_o} \\
=& \frac{2k^4(\om_o)}{\om_o^2 \gamma_{\om_o}}
\hspace{-0.05in}\left\{
\frac{1}{\sqrt{k^2(\om_o)-\mu^2(0)}} +\int_{z_T(\om_o)}^0 \hspace{-0.1in}dz \, \frac{\mu(z) [\mu(z) \mu'(z_T(\om_o))- \mu(z_T(\om_o)) \mu'(z)]}{k(\om_o)[k^2(\om_o)-\mu^2(z)]^{3/2}}
\right\}
\end{align*}
is an effective dispersion coefficient as we will see below.
When we observe $f_{\rm ref}^\ep$ around time $T_{\om_o}$, in a time window of order
$\sqrt{\ep}$, which corresponds to the scaled support of the emitted
pulse, we obtain
\begin{align}
f_{\rm ref}^\ep \left(T_{\om_o} + \sqrt{\ep} t\right) = \mbox{Re}\left\{ \exp
\left[ i \big(2 \phi_{\om_o}(0)-\om_oT_{\om_o} - \om_o \sqrt{\ep}
    t \big)/\ep\right] \mathcal{F}_{\rm ref}^\ep(t) \right\}.
\label{eq:PS4}
\end{align}
This oscillates at carrier frequency $\om_o/\ep$, like the emitted
pulse, and  its envelope
\begin{align}
\mathcal{F}_{\rm ref}^\ep(t)= \int \frac{d w}{2 \pi B} \hat F \Big( \frac{w}{B}\Big) \exp \left[
  i w^2 \beta_{\om_o} + i \psi_{\om_o + \sqrt{\ep} w}^\ep(0) - i w t \right],
\label{eq:PS5}
\end{align}
is described in the next theorem.

\vspace{0.05in}
\begin{theorem}
\label{thm.1}
Suppose that the standard deviation $\sigma_\ep$ of the random
fluctuations is of the order $|\ln \ep|^{-1/2}$, so that 
\begin{equation}
\upsilon_{\om_o}^2 =\frac{k^4(\om_o)}{\gamma_{\om_o}}
\hat \cR(0) \lim_{\ep \to 0} \sigma_\ep^2 \ln \left(
\frac{\phi_{\om_o}(0)}{\ep} \right) 
\label{eq:PS7}
\end{equation}
is finite, where $\gamma_{\om_o}$ is defined in \eqref{eq:res13} and
$\hat \cR(0) > 0$, because it is the power spectral density of the
fluctuations $\nu$, evaluated at zero. Then, as $\ep \to 0$,
$\mathcal{F}_{\rm ref}^\ep(t)$ converges in distribution, in the space of
continuous functions on compact sets in $\mathbb{R}$, to
\begin{equation}
\mathcal{F}_{\rm ref}(t) = \exp \Big( i \psi_{\om_o} - \frac{\upsilon_{\om_o}^2}{6}
\Big) \int \frac{d w}{2 \pi B} \, \hat F \Big( \frac{w}{B} \Big) \exp \big( i w^2
  \beta_{\om_o} - i w t \big),
\label{eq:PS6}
\end{equation}
where $\psi_{\om_o}$ is a Gaussian random variable with mean zero and
variance $2 \upsilon_{\om_o}^2/3$.
\end{theorem}

\vspace{0.05in} This is the pulse stabilization result, proved in
section \ref{sect:proof}. It says that aside from the random phase
$\psi_{\om_o} \sim \mathcal{N}(0,2 \upsilon_{\om_o}^2/3)$, the envelope of the
reflected pulse is deterministic.  It differs from the envelope $F(Bt)$
of the emitted pulse by the damping factor $\exp[-\upsilon_{\om_o}^2/6]$ and
the deformation by the second-order dispersive term $\beta_{\om_o} w^2$ in the phase.

\begin{remark}
\end{remark}The convergence stated in Theorem \ref{thm.1} holds in the space of continuous functions 
endowed with the topology induced by the supremum norm over the compact sets. It 
does not hold  in $L^2$.  Equation \eqref{eq:PS6} describes a reflected pulse with damped amplitude.
Its energy  is smaller than the energy of the incoming pulse, which may seem surprising because no energy 
can be transmitted beyond the turning point and there is no dissipation in the medium,
so all the incoming energy should be reflected. This is what we have before taking the limit $\ep \to 0$,
since the reflection coefficient has modulus one.
Theorem \ref{thm.1} describes only the coherent reflected pulse, which is observed around
the time $T_{\om_o}$, at the time scale of the incoming pulse width.
The theorem does not describe  the coda wave, consisting of the  incoherent, small-amplitude, long lasting  wave fluctuations that arrive after the coherent reflected pulse. These carry the remainder  of the energy.
When  $\upsilon_{\om_o} \ll 1$, the incoherent wave fluctuations are negligible, but when $\upsilon_{\om_o} \gg 1$, 
they carry most of the energy.

\begin{remark}
\end{remark}Theorem \ref{thm.1} assumes that the scaled standard
deviation $\sigma_\ep$ of the random fluctuations of the boundary is
small, of order $|\ln \ep|^{-1/2}$. In the absence of the turning
point such fluctuations would have a negligible effect on the
wave. The results \cite{alonso2011wave,gomez2011wave} in waveguides
with random perturbations of straight boundaries show that: (1) the
fluctuations have a net scattering effect when $\sigma_\ep = O(1)$,
and (2) the slower the modes propagate along the waveguide axis, the
stronger this effect.  We have a single propagating mode, which slows
down as it approaches the turning point, meaning that its group
velocity along $z$ tends to zero. Because the mode hovers around
$z_T(\om_o)$, it scatters repeatedly at the random boundary, which is
why the net scattering effect can be observed at the smaller standard
deviation $\sigma_\ep = O(|\ln \ep|^{-1/2})$.

Theorem \ref{thm.1} is proved in the following section.
Roughly speaking,
the proof is based on a diffusion-approximation result that describes the joint
distribution of the frequency-dependent phases of the 
reflection coefficients in the limit $\ep \to 0$.
These phases are random and become asymptotically Gaussian distributed,
and their covariance function (as a function of the frequency) exhibits an 
interesting feature: the frequency-dependent phases have a common random component
and they also have uncorrelated and identically distributed components.
Since the time-dependent profile of the reflected wave is the superposition of 
many frequency-dependent reflection coefficients by Fourier synthesis,
the common phase gives the random phase in the time-dependent profile
in (\ref{eq:PS6}), while the uncorrelated phases average out and give
 the damping term.

\section{Derivation of the pulse stabilization result}
\label{sect:proof}
We begin in section \ref{sect:proof1} with the single-frequency
asymptotic analysis of the random phase $\psi_\om^\ep(z)$, which
defines the reflection coefficient \eqref{eq:res35}. The
multi-frequency analysis of $\psi_\om^\ep(z)$ is in section
\ref{sect:proof2}, and the proof of Theorem \ref{thm.1} is completed
in section \ref{sect:proof3}. We assume throughout the section that 
$\sigma_\ep$  is of order $|\ln \ep|^{-1/2}$, as stated in
Theorem \ref{thm.1}.  
\subsection{Single-frequency analysis}
\label{sect:proof1}
To analyze the random phase $\psi_\om^\ep(z)$ in the limit $\ep \to
0$, we change variables so that equation \eqref{eq:res36ODE} takes a
form that can be analyzed with the diffusion limit theorem in
\cite{kim1996uniform}.  The change of variables is
\begin{equation}
z \to \zeta := \ep^{1/3} \eta_\om^\ep(z),
\label{eq:Pr1}
\end{equation}
with $\eta_\om^\ep$ defined in \eqref{eq:res8}. The inverse of this
mapping is $z =Z_\om (\ep^{1/3}\zeta)$ in terms of the function $Z_\om:\mathbb{R} \to
\mathbb{R}$, defined pointwise as the unique solution of
\begin{align}
\phi_\om \big(Z_\om(\xi)\big) = \frac{2}{3} \mbox{sgn}(\xi) \big| \xi
\big|^{3/2}, \qquad \forall \, \xi \in \mathbb{R}, \label{eq:Pr2}
\end{align}
with $\phi_\om$ given in \eqref{eq:res7}, and ``sgn'' denoting the
sign function.  Equivalently, in differential equation form,
$Z_\om(\xi)$ is the unique solution of
\begin{equation} 
\partial_\xi Z_\om(\xi) = Q^2_\om\big(Z_\om(\xi)\big) ~ ~\mbox{for}~~ \xi \ne
0, \qquad Z_\om(0) = z_T(\om).
\label{eq:Pr3}
\end{equation}

We denote the phase after the change of variables \eqref{eq:Pr1} with the same symbol $\psi_\om^\ep$, 
and obtain using  \eqref{eq:res14} that equation
\eqref{eq:res36ODE} becomes
\begin{align}
\partial_\zeta \psi_\om^\ep(\zeta) &=\frac{2J_\om^2(\ep^{1/3} \zeta) V(\ep^{-1/3} \zeta)}{
  \sqrt{|\zeta|}} 
\left[ \frac{\sigma_\ep}{\ep^{1/3}} \mu^2 \big(Z_\om(\ep^{1/3} \zeta)\big) \nu \Big(
\frac{Z_\om(\ep^{1/3} \zeta)}{\ep} \Big)  \right. \nonumber \\
&\hspace{-0.3in}\left. +\sigma_\ep^2 \ep^{1/6} g^\ep \big(\om, Z_\om(\ep^{1/3} \zeta)  \big) \right]  \cos^2
\Big\{ \frac{\psi_\om^\ep(\zeta)}{2} - \arg[A_i + i B_i](-\ep^{-1/3} \zeta)\Big\},
\label{eq:Pr4LONG}
\end{align}
with 
\begin{equation}
J_\om(\xi) = Q_\om^2(Z_\om(\xi)), \qquad V(\xi) = \pi \sqrt{|\xi|} \left[
  A_i^2(-\xi) + B_i^2(-\xi)\right],
\label{eq:Pr5}
\end{equation}
and the shortened notation 
\[
\arg[A_i+i B_i](-\xi) = \arg[A_i(-\xi) + i B_i(-\xi)].
\]
The turning point lies at $\zeta = 0$, and the source is at
\begin{equation}
\label{eq:zeta0}
\zeta_{\rm s}^\ep = \ep^{1/3} \eta_\om^\ep(0) = \ep^{-1/3} [3
  \phi_\om(0)/2]^{2/3} = O(\ep^{-1/3}).
\end{equation}
The boundary point $z_{\rm b}$ is mapped to
\begin{equation}
\label{eq:zeta_b}
\zeta_{\rm b}^\ep = -\ep^{-1/3} [-3 \phi_\om(z_{\rm b})/2]^{2/3},
\end{equation}
and we have the boundary condition 
\begin{equation}
\psi_\om^\ep(\zeta_{\rm b}^\ep) = 0.
\label{eq:Pr6}
\end{equation}
The quadratic term in the fluctuations, modeled by $g^\ep$ in \eqref{eq:Pr4LONG}, is 
negligible in the limit $\ep \to 0$, 
because it gives a contribution that can be bounded by 
\[O( \sigma_\ep^2 \ep^{1/6}\sqrt{\zeta^\ep_{\rm s}-\zeta^\ep_{\rm b}} )=
O(\sigma_\ep^2) = O(1/|\ln \ep|).
\]
We neglect it henceforth and simplify equation \eqref{eq:Pr4LONG} to 
\begin{align}
\partial_\zeta \psi_\om^\ep(\zeta) =& \frac{2 \sigma_\ep}{\ep^{1/3}
  \sqrt{|\zeta|}} J_\om^2(\ep^{1/3} \zeta) V(\ep^{-1/3} \zeta)
\mu^2 \big(Z_\om(\ep^{1/3} \zeta)\big) \nu \left(
\frac{Z_\om(\ep^{1/3} \zeta)}{\ep} \right) \nonumber \\ &\times \cos^2
\Big\{ \frac{\psi_\om^\ep(\zeta)}{2} - \arg[A_i+ i B_i](-\ep^{-1/3} \zeta)\Big\}. 
\label{eq:Pr4}
\end{align}

To understand how $\psi_\om^\ep(\zeta)$ evolves from the
boundary value \eqref{eq:Pr6}, let us start from $\zeta_{\rm b}^\ep$ and
consider first points that are far on the left of the turning point, 
at $\zeta < 0$ satisfying $|\zeta| \gg O(\ep^{1/3})$. The function $V$
defined in \eqref{eq:Pr5} is large at these points, as given by the
asymptotic expansions of the Airy function $B_i$ in appendix
\ref{ap:propagator.3},
\begin{equation*}
V(\ep^{-1/3} \zeta) \approx \pi \sqrt{\ep^{-1/3} |\zeta|} \,
B_i^2(\ep^{-1/3}|\zeta|) \approx \exp \left[\frac{4
    |\zeta|^{3/2}}{3 \ep^{1/2}}\right].
\end{equation*}
We also have the expansion 
\begin{equation*}
\arg\left[A_i(-\ep^{-1/3} \zeta) + i B_i(-\ep^{-1/3} \zeta)\right]
\approx \frac{\pi}{2} - \frac{1}{2} \exp \left[-\frac{4
    |\zeta|^{3/2}}{3 \ep^{1/2}}\right],
\end{equation*}
and since the phase starts from zero by \eqref{eq:Pr6}, 
\begin{align*}
\cos^2\Big\{ \frac{\psi_\om^\ep(\zeta)}{2} - \arg\left[A_i(-\ep^{-1/3}
  \zeta) + i B_i(-\ep^{-1/3} \zeta)\right]\Big\} \approx 
\frac{1}{4} \exp \left[-\frac{8 |\zeta|^{3/2}}{3 \ep^{1/2}}\right].
\end{align*}
This makes the right-hand side in \eqref{eq:Pr4} exponentially small,
so the phase remains essentially zero on the left of the turning
point. Moreover, the phase is independent of the precise value of
$\zeta_{\rm b}^\ep$  at which we prescribe the boundary condition
\eqref{eq:Pr6}.

Now consider the $O(\ep^{1/3})$ vicinity of the turning point, where
we can set $\zeta = \ep^{1/3} \tilde{\zeta}$ with $\tilde{\zeta} =
O(1)$, to rewrite equation \eqref{eq:Pr4} for $\widetilde
\psi_\om^\ep(\tilde{\zeta}) = \psi_\om^\ep(\ep^{1/3}\tilde{\zeta})$ as
\begin{align}
\partial_{\tilde{\zeta}} \widetilde \psi_\om^\ep(\tilde{\zeta}) =
\frac{2 \pi \sigma_\ep k^2(\om)}{\ep^{1/6} \gamma_\om^{2/3}}
\left[A_i^2(-\tilde{\zeta}) + B_i^2(-\tilde{\zeta})\right] \nu \left(
\frac{z_T(\om)}{\ep} + \frac{\tilde{\zeta}}{\ep^{1/3}
  \gamma_\om^{1/3}} \right) \nonumber \\ \times \cos^2 \left\{
\frac{\widetilde \psi_\om^\ep(\tilde{\zeta})}{2} - \arg
[A_i + i B_i](-\tilde{\zeta}) \right\} +
\ldots,
\label{eq:Pr7p}
\end{align}
with the dots denoting negligible terms. Here we used equation
\eqref{eq:Pr3}, and 
\[J_\om(0) = \gamma_\om^{-1/3}, \qquad \mu^2(Z_\om(0)) =
\mu^2(z_T(\om)) = k^2(\om).
\]  
If $\sigma_\ep$ were order one, the right-hand side in \eqref{eq:Pr7p}
would be in the usual diffusion approximation form with 
\[
\widetilde \nu^\ep(\tilde \zeta) =  \frac{1}{\ep^{1/6}}
\nu\left(\frac{z_T(\om)}{\ep}+ \frac{\tilde{\zeta}}{\ep^{1/3}
  \gamma_\om^{1/3}} \right) 
\] behaving like white noise in the limit
$\ep \to 0$, for $\tilde{\zeta}$ of order one
\cite{papanicolaou1974asymptotic}. But in our case $\sigma_\ep =
O(|\ln \ep|^{-1/2})$ tends to zero as $\ep \to 0$, so the fluctuations
are negligible in the $O(\ep^{1/3})$ vicinity of the turning point.

The net scattering effect at the random boundary comes from the long
interval
\begin{equation*}
\mathcal{I}^\ep = \left\{ \zeta \in \mathbb{R}~~ \mbox{s.t.} ~ ~ O(\ep^{1/3}) < 
\zeta \le \zeta_{\rm s}^\ep \right\},
\end{equation*} 
that grows as $\ep^{-1/3}$ in the limit $\ep \to 0$ by (\ref{eq:zeta0}). The asymptotic
expansions of the Airy functions at large negative arguments given in
appendix \ref{ap:propagator.2} show that in $\mathcal{I}^\ep$ we have 
\begin{equation}
V(\ep^{-1/3} \zeta) \approx 1 + o(1),
\label{eq:Pr7}
\end{equation}
and 
\begin{equation}
\arg\left[A_i(-\ep^{-1/3} \zeta) + i B_i(-\ep^{-1/3} \zeta)\right]
\approx \frac{\pi}{4} - \frac{2}{3}(\ep^{-1/3} \zeta)^{3/2} + o(1).
\label{eq:Pr8}
\end{equation}
One can verify that these are excellent approximations for
all $\ep^{-1/3} \zeta >3$, so we can take $\mathcal{I}^\ep =
(\zeta_{-}^\ep,\zeta_{\rm s}^\ep]$, with $\zeta_{-}^\ep = 3 \ep^{1/3}$, and
  simplify equation \eqref{eq:Pr4} as
\begin{align}
\partial_\zeta \psi_\om^\ep(\zeta) = \frac{\sigma_\ep
  J_\om^2(\ep^{1/3} \zeta)\mu^2 \big(Z_\om(\ep^{1/3}
  \zeta)\big)}{\ep^{1/3} \sqrt{\zeta}} \nu \left(
\frac{Z_\om(\ep^{1/3} \zeta)}{\ep} \right) \nonumber \\ \times
\Big\{1 + \sin\Big[ \psi_\om^\ep(\zeta) + \frac{4}{3}(\ep^{-1/3}
  \zeta)^{3/2}\Big]\Big\}.
\label{eq:Pr9}
\end{align}

\subsection{The diffusion limit for the single-frequency case}
\label{sect:diflim1}
To see how to apply the limit theorem in \cite{kim1996uniform} to
equation \eqref{eq:Pr4}, imagine that we discretize the interval
$\mathcal{I}^\ep$ at points $\zeta^{(j)} = \zeta_{-}^\ep + j
\Delta_\zeta$ separated by $\Delta_\zeta = O(1)$, for $j = 0, \ldots,
n^\ep$, and 
\begin{equation} 
\label{eq:Pr10a}
n^\ep = \lfloor(\zeta_{\rm s}^\ep-\zeta_{-}^\ep)/\Delta_\zeta\rfloor = O(\ep^{-1/3}). 
\end{equation}
If the argument of $\lfloor ~\rfloor$ in this equation is not integer,
the length of the last interval is adjusted so that 
$\zeta_{\rm s}^\ep = \zeta^{(n^\ep)}$.

In each sub-interval $[\zeta^{(j)},\zeta^{(j+1)}]$, we
expand the argument of $\nu$ in \eqref{eq:Pr9} as
\begin{equation*}
\frac{Z_\om(\ep^{1/3} \zeta)}{\ep} \approx \frac{Z_\om(\ep^{1/3}
  \zeta^{(j)})}{\ep} + J_\om\big(\ep^{1/3} \zeta^{(j)}\big)
\frac{\zeta-\zeta^{(j)}}{\ep^{2/3}} +
\left(\frac{\zeta-\zeta^{(j)}}{\ep^{2/3}}\right)^2 \hspace{-0.05in}O(\ep),
\end{equation*}
where $\ep^{1/3} \zeta^{(j)}$ is order one by \eqref{eq:Pr10a}, and we
used equations \eqref{eq:Pr3} and \eqref{eq:Pr5}. Similarly, the
argument of the sin in \eqref{eq:Pr9} is
\begin{equation*}
(\ep^{-1/3} \zeta)^{3/2} \approx \big(\ep^{-1/3} \zeta^{(j)}\big)^{3/2} +
  \frac{3 \sqrt{\ep^{1/3} \zeta^{(j)}}}{2}
  \frac{\zeta-\zeta^{(j)}}{\ep^{2/3}} + \left(
  \frac{\zeta-\zeta^{(j)}}{\ep^{2/3}}\right)^2\hspace{-0.05in}
\frac{O(\ep)}{\sqrt{\ep^{1/3}
      \zeta^{(j)}}}.
\end{equation*}
This makes equation \eqref{eq:Pr9} of the same form as in
\cite{kim1996uniform}, with the small parameter $\epsilon$ there
replaced by our $\ep^{1/3}$, and the process $\nu$ satisfying by
assumption the strong mixing conditions in \cite{kim1996uniform}.
Thus, we can use the limit theorem in \cite[section
  III]{kim1996uniform} for the joint process $\left(\nu^\ep(\zeta),
\zeta^\ep(\zeta) \right)$ on the state space $\mathbb{R} \times [0,3
  \pi/2]$, with
\begin{equation}
\label{eq:Pr10b}
\nu^\ep(\zeta) = \nu \left(\frac{Z_\om(\ep^{1/3} \zeta)}{\ep}\right),\qquad 
\zeta^\ep(\zeta) = \left(\ep^{-1/3}\zeta \right)^{3/2}.
\end{equation}
The torus $[0,3 \pi/2]$ arises because the right-hand side in \eqref{eq:Pr9} is 
periodic in $\zeta^\ep(z)$.

The next lemma, proved in appendix \ref{ap:diflimit},  describes the distribution of the phase $\psi_\om^\ep(\zeta_{\rm s}^\ep)$ at the 
location $\zeta_{\rm s}^\ep$ of the source. 
\vspace{0.05in}
\begin{lemma}
\label{lem.D1}
$\psi_\om^\ep(\zeta_{\rm s}^\ep)$ is asymptotically Gaussian distributed in
the limit $\ep \to 0$, with mean zero and variance 
\begin{equation}
\upsilon_\om^2 = \frac{k^4(\om)}{\gamma_\om} \hat \cR(0) \lim_{\ep \to
  0} \sigma_\ep^2 \ln \left(\frac{\phi_\om(0)}{\ep}\right).
\label{eq:Pr11}
\end{equation}
\end{lemma}

\subsection{Multi-frequency analysis}
\label{sect:proof2}
The expression \eqref{eq:PS5} of the envelope of the reflected pulse
involves the random phases $\psi_{\om}^\ep(z=0)$ at frequencies $\om =
\om_o + \sqrt{\ep} w$, with $w$ in the support of $\hat F(w/B)$, the
Fourier transform of the envelope of the emitted pulse. Here we
describe the asymptotic distribution of these phases at $m$
such distinct frequencies.

Let us introduce the notation
\begin{equation}
\psi_j^\ep(z) = \psi_{\om_o + \sqrt{\ep} w_j}^\ep(z),
\label{eq:Pr12}
\end{equation}
and consider the random process $\left(\psi_1^\ep(z), \ldots ,\psi_m^\ep(z)\right)$. 
Each $\psi_j^\ep(z)$ satisfies equation
\eqref{eq:res36ODE} with $\om = \om_o + \sqrt{\ep} w_j$, and boundary
condition $\psi_j^\ep(z_{\rm b}) = 0$. We proceed as in the previous section
and change variables to transform the problem into one that can be
analyzed with the diffusion limit theorem in
\cite{kim1996uniform}. The change of variables is similar to
\eqref{eq:Pr1}, but since we have multiple frequencies that are close
to $\om_o$, we take
\begin{equation}
z \to \zeta := \ep^{1/3} \eta_{\om_o}^\ep(z),
\label{eq:Pr13}
\end{equation}
with inverse transform given by $z=Z_{\om_o}(\ep^{1/3}\zeta)$ with $Z_{\om_o}$ defined in \eqref{eq:Pr2}
for $\om = \om_o$. By definitions \eqref{eq:res8} and \eqref{eq:Pr2} we 
have the expansion
\begin{align}
\eta_{\om_o + \sqrt{\ep} w}^\ep(Z_{\om_o}(\xi)) &=
\mbox{sgn}(\xi)\ep^{-2/3}\left[ \mbox{sgn}(\xi) \frac{3}{2}
  \phi_{\om_o + \sqrt{\ep}w}(Z_{\om_o}(\xi)) \right]^{2/3} \nonumber
\\ &= \ep^{-2/3} \left[ \xi + \ep^{1/2} w \mathcal{K}(\xi) \right] +
o(1), \label{eq:Pr14}
\end{align}
uniformly in $\xi$, up to $|\xi| = O(1)$, where
\begin{equation}
\mathcal{K}(\xi) = 
\frac{1}{|\xi|^{1/2}} \partial_{\om}
\phi_{\om}(Z_{\om_o}(\xi))\big|_{\om = \om_o} .
\end{equation}
The function 
\begin{equation} 
\mathcal{K}(\xi) = 
\left\{
\begin{array}{rr}
\displaystyle
\frac{k^2(\om_o)}{\om_o \xi^{1/2}}
\int_{z_T(\om_o)}^{Z_{\om_o}(\xi)} \frac{d z}{\sqrt{k^2(\om_o) -
    \mu^2(z)}}, & \mbox{ if $\xi >0$,} \\
\displaystyle
\frac{k^2(\om_o)}{\om_o |\xi|^{1/2}}
\int_{Z_{\om_o}(\xi)}^{z_T(\om_o)}  \frac{d z}{\sqrt{\mu^2(z)-k^2(\om_o) }}, & \mbox{ if $\xi<0$},
    \end{array}
    \right.
\end{equation}
is continuous and it is equal to $2 k^2(\om_o)/(
\gamma_{\om_o}^{2/3}\om_o)$ at the turning point, where $\xi = 0$.

Using the change of variables \eqref{eq:Pr13} and the expansion
\eqref{eq:Pr14} for $\xi = \ep^{1/3} \zeta$ in \eqref{eq:res36ODE} we
obtain
\begin{align}
\partial_\zeta \psi_j^\ep(\zeta) =& \frac{2 \sigma_\ep}{\ep^{1/3}
  \sqrt{|\zeta|}} J_{\om_o}^2(\ep^{1/3} \zeta) V(\ep^{-1/3} \zeta)
\mu^2 \big(Z_{\om_o}(\ep^{1/3} \zeta)\big) \nu \left(
\frac{Z_{\om_o}(\ep^{1/3} \zeta)}{\ep} \right) \nonumber \\ &\times
\cos^2 \Big\{ \frac{\psi_\om^\ep(\zeta)}{2} - \arg\left[A_i + i B_i
  \right]\big(-\ep^{-1/3} \zeta - \ep^{-1/6} w_j
\mathcal{K}(\ep^{1/3} \zeta))\Big\},
\label{eq:Pr15}
\end{align}
where we neglected the terms that have no contribution as $\ep \to 0$. 
This equation  is basically the same as equation
\eqref{eq:Pr4} analyzed in the previous section, except that the
phases in the argument of the cosine change with $j$.

The discussion in the previous section applies verbatim here, and we
conclude the same way that we need to consider only $\zeta \in
(\zeta_{-}^{\ep}, \zeta_{\rm s}^\ep]$, with 
\[
\zeta_{-}^\ep = 3 \ep^{1/3}, \qquad \zeta_{\rm s}^\ep = \ep^{-1/3} \left[3
  \phi_{\om_o}(0)/2\right]^{2/3} = O(\ep^{-1/3}),
\] 
where
\eqref{eq:Pr15} takes the form
\begin{align}
\partial_\zeta \psi_j^\ep(\zeta) =& \frac{\sigma_\ep}{\ep^{1/3}
  \sqrt{|\zeta|}} J_{\om_o}^2(\ep^{1/3} \zeta) \mu^2
\big(Z_{\om_o}(\ep^{1/3} \zeta)\big) \nu \left(
\frac{Z_{\om_o}(\ep^{1/3} \zeta)}{\ep} \right) \nonumber \\ &\times
\left\{ 1 + \sin \left[ \psi_j^\ep(\zeta) + \frac{4}{3}\big(\ep^{-1/3}
  \zeta + \ep^{-1/6} w_j \mathcal{K}(\ep^{1/3}
  \zeta)\big)^{3/2}\right]\right\}.
\label{eq:Pr16}
\end{align}
This leads to the asymptotic distribution of the phases stated in the next lemma,
proved in  appendix \ref{ap:diflimit2}.

\vspace{0.05in}
\begin{lemma} 
\label{lem.D2} The vector 
\begin{equation}
\boldsymbol{\Psi}^\ep(\zeta_{\rm s}^\ep) = \left(\psi_1^\ep(\zeta_{\rm s}^\ep), 
\ldots, \psi_m^\ep(\zeta_{\rm s}^\ep) \right), \label{eq:process}
\end{equation}
converges in distribution in the limit $\ep \to 0$ to a Gaussian vector with mean zero
and covariance matrix
\begin{equation}
{\bf C} = \frac{\upsilon_{\om_o}^2}{3} ({\bf I}_m + 2 {\bf J}_m),
\label{eq:covar}
\end{equation}
where $\upsilon_{\om_o}^2$ is defined in \eqref{eq:PS7}, ${\bf I}_m$ is the $m \times m$
identity matrix, and ${\bf J}_m$ is the $m \times m$
matrix with all entries equal to one.
\end{lemma}
\subsection{Proof of pulse stabilization}
\label{sect:proof3}

By Lemma \ref{lem.D2} and  definition
\eqref{eq:PS5} of the envelope of the reflected pulse, 
for any $t_1,\ldots,t_m \in \RR$ we can  
calculate the finite-order moments
\begin{align}
\EE \left[ \prod_{j=1}^m \mathcal{F}_{\rm ref}^\ep(t_j) \right] = &\int
\frac{d w_1}{2 \pi B} \hat F \Big( \frac{w_1}{B} \Big) \ldots \int \frac{d w_m}{2 \pi B} \hat
F\Big( \frac{w_m}{B} \Big) \exp \Big[i \sum_{j=1}^m \left( w_j^2 \beta_{\om_o} - w_j
  t_j\right) \Big]\nonumber \\ &\times \EE \Big[\exp \Big(i \sum_{j=1}^m
  \psi_j^\ep(\zeta_{\rm s}^\ep) \Big) \Big],
\label{eq:Pf1}
\end{align}
in the limit $\ep \to 0$. The process \eqref{eq:process} is Gaussian
in this limit, with covariance \eqref{eq:covar}, so the expectation in
\eqref{eq:Pf1} is given by
\begin{equation}
\lim_{\ep \to 0} \EE \Big[\exp \Big(i \sum_{j=1}^m
  \psi_j^\ep(\zeta_{\rm s}^\ep) \Big) \Big] = \exp \left[ -\frac{m(2m+1)
    \upsilon_{\om_o}^2}{6} \right].
\label{eq:Pf2}
\end{equation}
The right-hand side can also be written in terms of the random Gaussian phase $\psi_{\om_o}$
with mean zero and variance $2 \upsilon_{\om_o}^2/3$ as 
\begin{equation}
\exp \left[ -\frac{m(2m+1) \upsilon_{\om_o}^2}{6} \right] = \EE\Big[ \prod_{j=1}^m
\exp \Big(-\frac{\upsilon_{\om_o}^2}{6} + i 
  \psi_{\om_o}\Big)\Big] .
\label{eq:Pf3}
\end{equation}
By substituting into \eqref{eq:Pf1}, we obtain the convergence of the finite-order moments
\begin{equation}
\lim_{\ep \to 0} \EE \left[ \prod_{j=1}^m \mathcal{F}_{\rm ref}^\ep(t_j)
  \right] = \EE \left[ \prod_{j=1}^m \mathcal{F}_{\rm ref}(t_j) \right],
\label{eq:Pf4}
\end{equation}
for $\mathcal{F}_{\rm ref}(t)$ defined in \eqref{eq:PS6}. 

The convergence result stated in Theorem \ref{thm.1} follows from
\eqref{eq:Pf4}, once we prove tightness of the process
$\mathcal{F}_{\rm ref}^\ep(t)$ in the
space of continuous functions on compact sets in $\mathbb{R}$, as
shown in \cite[Chapter 2]{billingsley2013convergence}. 

We obtain from definition \eqref{eq:PS5} and the triangle inequality
that $F_{\rm ref}^\ep(t)$ is bounded independent of $\ep$, uniformly in
$t \in [0,T]$, where $T$ is finite. Moreover, for any $\Delta_t > 0$, 
\begin{align}
\left|\mathcal{F}_{\rm ref}^\ep(t + \Delta_t)-\mathcal{F}_{\rm ref}^\ep(t)
\right|&= \Big|\int \frac{dw}{2 \pi B} \, \hat F \Big( \frac{w}{B} \Big) e^{ i w^2 \beta_{\om_o} +
  i \psi_{\om_o + \sqrt{\ep} w}^\ep(\zeta_{\rm s}^\ep) - i w (t+\Delta_t)} \left(1 -
e^{i w \Delta_t}\right)\Big| \nonumber \\ &\le \int \frac{dw}{2 \pi B} \,
\Big|\hat F \Big(\frac{w}{B} \Big) \Big| \big|1 - \exp(i w \Delta_t)\big|.
\label{eq:Pf5}
\end{align}
Note that 
\[
\big|1 - e^{i w \Delta_t}\big| = 2 |\sin (w \Delta_t / 2)| \le  |w| \Delta_t,
\]
and that $w \hat F(w/B)$ is absolutely integrable by the assumption of
compact support $[-\pi,\pi]$ of $\hat F(w)$. Substituting in
\eqref{eq:Pf5}, we conclude that there exists a constant $C$,
independent of $\ep$ and $t$, that bounds the modulus of continuity
\begin{align*}
\sup_{|t'-t| \le
  \Delta_t} \left|\mathcal{F}_{\rm ref}^\ep(t')-\mathcal{F}_{\rm ref}^\ep(t)
\right| \le C \Delta_t.
\end{align*}
This implies that $\big(\mathcal{F}_{\rm ref}^\ep(t)\big)_{t \in
  \mathbb{R}}$ is tight 
   \cite[Chapter 2]{billingsley2013convergence}.
This completes the proof of Theorem \ref{thm.1}. $\Box$.
\section{Summary}
\label{sect:sum}
In this paper we studied the reflection of a pulse in a random waveguide with a turning point. 
The waveguide has reflecting boundaries, a slowly bending axis, and variable cross section. 
The variation consists of small-amplitude, random fluctuations
of the boundary, and a slow and monotone change of the  opening of the waveguide. 
The pulse is emitted by a point source, and is modeled as usual by a carrier oscillatory signal 
multiplying a smooth envelope. The carrier wavelength is similar to the 
width of the cross section of the  waveguide, so that the emitted wave is a superposition of a single propagating mode and 
 infinitely many evanescent modes.  The  
turning point is many wavelengths away from the source, and 
marks the limit of propagation of the mode in the waveguide, meaning that once the wave reaches it, it is reflected back. 
The goal of the paper is to characterize in detail  this reflection. 

We derived from first principles, starting with the wave equation in the waveguide, a stochastic differential equation 
for the reflection coefficient, driven by the random fluctuations of the boundary. We showed how
this equation can be studied asymptotically, for a small carrier wavelength with respect to the distance of propagation,
using stochastic diffusion limits. We also quantified the amplitude of the random fluctuations of the 
boundary under which the reflected pulse is strongly affected and we explain why it
maintains a deterministic shape i.e., it is stabilized. The reflected
pulse oscillates at the same central frequency as the emitted one, but it has a different envelope that is 
damped and deformed due to scattering at the random boundary.

\section*{Acknowledgements}
Liliana Borcea's work was partially supported by the NSF grant
DMS1510429. Support from AFOSR grant FA9550-15-1-0118 is also 
gratefully acknowledged.

\appendix
\section{Useful identities}
\label{ap:identities}
Here we give a few identities satisfied by the eigenfunctions
\eqref{eq:eigf}, for all $z \in \mathbb{R}$. The first identity
is just the statement that the eigenfunctions are orthonormal 
\begin{equation}
\int_{-D(z)/2}^{D(z)/2} d \rho \, y_j(\rho,z) y_q(\rho,z) =
\delta_{jq},
\label{eq:orthog}
\end{equation}
where $\delta_{jq}$ is the Kronecker delta symbol. The second identity
\begin{align}
\int_{-D(z)/2}^{D(z)/2} d \rho \, \rho y_j^2(\rho,z)&= 0,
\label{eq:odd}
\end{align}
is due to the fact that the integrand is odd. The third identity
follows from the fundamental theorem of calculus,
\begin{align}
\int_{-D(z)/2}^{D(z)/2} d \rho\,y_j(\rho,z) \partial_\rho
y_j(\rho,z) &= \frac{1}{2}\int_{-D(z)/2}^{D(z)/2} d \rho \,
\partial_\rho y_j^2(\rho,z) = 0,
\end{align}
because the eigenfunctions vanish at $\rho = \pm D(z)/2$. 
The fourth identity is 
\begin{align}
\int_{-D(z)/2}^{D(z)/2} d \rho \, [2\rho + D(z)] y_j(\rho,z)
\partial_\rho y_j(\rho,z) = \int_{-D(z)/2}^{D(z)/2} d \rho \, \rho
\partial_\rho y_j^2(\rho,z) \nonumber \\ =\int_{-D(z)/2}^{D(z)/2} d
\rho \left\{ \partial_\rho \left[\rho y_j^2(\rho,z)\right] -
y_j^2(\rho,z)\right\} = -1,
\end{align}
where we used integration by parts. The fifth identity is 
\begin{equation}
\int_{-D(z)/2}^{D(z)/2} d \rho \, y_j(\rho,z) \partial_z y_j(\rho,z) = 0,
\end{equation}
and to derive it, we take the derivative with respect to $z$ in
\eqref{eq:orthog}, for $q=j$, and obtain that
\begin{align*}
0 =& \partial_z \int_{-D(z)/2}^{D(z)/2} d \rho \, y_j^2(\rho,z) = 2
\int_{-D(z)/2}^{D(z)/2} d \rho \, y_j(\rho,z) \partial_z
y_j(\rho,z) \\&+ \frac{D'(z)}{2} \left[ y_j^2(D(z)/2,z)-
  y_j^2(-D(z)/2,z)\right] = 2
\int_{-D(z)/2}^{D(z)/2} d \rho \, y_j(\rho,z) \partial_z
y_j(\rho,z).
\end{align*}
The last identity follows from \eqref{eq:orthog}, \eqref{eq:odd}, and
the substitution of \eqref{eq:eigf} in the remaining integral that is
evaluated explicitly
\begin{align}
\int_{-D(z)/2}^{D(z)/2} d \rho [2\rho + D(z)]^2 y_j^2(\rho,z) &=
D^2(z) + \frac{8}{D(z)} \int_{-D(z)/2}^{D(z)/2} d \rho\, \rho^2 \sin^2
\left[ \left(\frac{\rho}{D(z)} + \frac{1}{2} \right) \pi j\right]
\nonumber \\&= D^2(z) \left[\frac{4}{3}-\frac{2}{(\pi j)^2}\right].
\end{align}
\section{Properties of the propagator}
\label{ap:propagator} 
Here we prove the statements of Lemmas \ref{lem.1}--\ref{lem.3}, which
describe the approximate propagator $\bM^\ep(\om,z)$.
\subsection{Proof of Lemma \ref{lem.1}} To derive equation
\eqref{eq:res17}, we need to show that
\begin{equation}
\partial_z M_{11}^\ep(\om,z) = \frac{i}{\ep} M_{21}^\ep(\om,z), \label{eq:prop1}
\end{equation}
and 
\begin{equation}
\partial_z M_{21}^\ep(\om,z) = \frac{i}{\ep} [k^2(\om)-\mu^2(z)] M_{11}^\ep(\om,z) -
i \ep \frac{Q_\om''(z)}{Q_\om(z)} M_{11}^\ep(\om,z),
\label{eq:prop2}
\end{equation}
for all $z < 0$. Equation \eqref{eq:prop1} is just the definition
\eqref{eq:res16} of $M_{21}^\ep(\om,z)$. Equation \eqref{eq:prop2}
follows by taking the derivative in the right-hand side of
\eqref{eq:res16}, using \eqref{eq:res14} and the equations satisfied
by the Airy functions \cite[chapter 10]{abramowitz1972handbook}
\begin{equation}
A_i''(-\eta_\om^\ep(z)) = - \eta_\om^\ep(z) A_i(-\eta_\om^\ep(z)), \qquad 
B_i''(-\eta_\om^\ep(z)) = - \eta_\om^\ep(z) B_i(-\eta_\om^\ep(z)).
\label{eq:prop3}
\end{equation}

The determinant of the matrix $\bM^\ep(\om,z)$ defined in \eqref{eq:res6} is given by
\begin{equation}
\det \bM^\ep(\om,z) = 2 \, \mbox{Re} \left[ M_{11}^\ep(\om,z)
  \overline{M_{21}^\ep(\om,z)} \right],
\label{eq:prop4}
\end{equation}
where $\mbox{Re}[~ ]$ denotes the real part. We calculate using
definitions \eqref{eq:res15}--\eqref{eq:res16} that
\begin{align*}
 M_{11}^\ep(\om,z) \overline{M_{21}^\ep(\om,z)} =& - i \pi
 \left[A_i(-\eta_\om^\ep(z))-iB_i(-\eta_\om^\ep(z))\right]
 \left[A_i'(-\eta_\om^\ep(z))+iB_i'(-\eta_\om^\ep(z))\right] \nonumber
 \\ &+ i \pi \ep^{2/3} Q_\om(z) Q'_\om(z) \left[A_i^2(-\eta_\om^\ep(z))
   + B^2_i(-\eta_\om^\ep(z))\right],
\end{align*}
and taking the real part we have 
\begin{align*}
\mbox{Re} \left[ M_{11}^\ep(\om,z) \overline{M_{21}^\ep(\om,z)}
  \right] = \pi \left[A_i(-\eta_\om^\ep(z)) B_i'(-\eta_\om^\ep(z)) - 
A_i'(-\eta_\om^\ep(z)) B_i(-\eta_\om^\ep(z))\right].
\end{align*}
The term in the square bracket in the right-hand side is the Wronskian of
the Airy functions which is constant and equal to $1/\pi$. Equation
\eqref{eq:res18} follows after substituting the result in \eqref{eq:prop4}. 
$\Box$

\subsection{Proof of Lemma \ref{lem.2}}
\label{ap:propagator.2} When $z \nearrow 0$, the 
function \eqref{eq:res7} is order one and $\eta_\om^\ep(z)$ defined in 
\eqref{eq:res8} satisfies
$
\eta_\om^\ep(z) = O( \ep^{-2/3}) \gg 1.
$
The asymptotic expansions of the Airy functions at
large and negative argument  \cite[chapter
  10]{abramowitz1972handbook} give that 
\begin{align}
A_i(-\eta_\om^\ep(z)) &\approx \frac{1}{\sqrt{\pi}
  \big(\eta_\om^\ep(z)\big)^{1/4}} \left\{\sin \left[ \frac{2}{3}
  \big(\eta_\om^\ep(z)\big)^{3/2} + \frac{\pi}{4} \right]+
  O\left(\big(\eta_\om^\ep(z)\big)^{-3/2}\right)\right\} ,
\\ A_i'(-\eta_\om^\ep(z)) &\approx
-\frac{\big(\eta_\om^\ep(z)\big)^{1/4}}{\sqrt{\pi}} \left\{ \cos
\left[ \frac{2}{3} \big(\eta_\om^\ep(z)\big)^{3/2} + \frac{\pi}{4}\right] +
  O\left(\big(\eta_\om^\ep(z)\big)^{-3/2}\right)\right\},
\end{align}
and 
\begin{align}
B_i(-\eta_\om^\ep(z)) &\approx \frac{1}{\sqrt{\pi}
  \big(\eta_\om^\ep(z)\big)^{1/4}}\left\{ \cos \left[ \frac{2}{3}
  \big(\eta_\om^\ep(z)\big)^{3/2} + \frac{\pi}{4} \right] +
O\left(\big(\eta_\om^\ep(z)\big)^{-3/2}\right) \right\},
\\ B_i'(-\eta_\om^\ep(z)) &\approx
\frac{\big(\eta_\om^\ep(z)\big)^{1/4}}{\sqrt{\pi}}\left\{ \sin \left[
  \frac{2}{3} \big(\eta_\om^\ep(z)\big)^{3/2} + \frac{\pi}{4}\right]+
  O\left(\big(\eta_\om^\ep(z)\big)^{-3/2}\right)\right\}.
\end{align}
We also get from definitions \eqref{eq:res8} and \eqref{eq:res9} that 
\begin{equation}
\big(\eta_\om^\ep(z)\big)^{1/4} = \ep^{-1/6} [k^2(\om)-\mu^2(z)]^{1/4} Q_\om(z),
\end{equation}
and 
\begin{equation}
\frac{2}{3} \big(\eta_\om^\ep(z)\big)^{3/2} = \ep^{-1} \phi_\om(z).
\end{equation}
The statement of the lemma follows by straightforward calculations
from these results and definitions \eqref{eq:res15}--\eqref{eq:res16}. $\Box$

\subsection{Proof of Lemma \ref{lem.3}}
\label{ap:propagator.3} When $z< z_T(\om)$ and 
$|z_T(\om)-z| \gg O(\ep^{2/3})$, the function $\eta_\om^\ep(z)$
defined in \eqref{eq:res8} is negative valued and 
$
|\eta_\om^\ep(z)|  \gg 1.
$
Then, we have from the asymptotic expansions of the Airy functions at
large and positive argument that \cite[chapter
  10]{abramowitz1972handbook}
\begin{align}
A_i(|\eta_\om^\ep(z)|) &\approx \frac{1}{2\sqrt{\pi}
  |\eta_\om^\ep(z)|^{1/4}} e^{-\frac{2}{3}
  \big|\eta_\om^\ep(z)\big|^{3/2}} \left[1 +
  O\left(\big|\eta_\om^\ep(z)\big|^{-3/2}\right)\right] ,
\\ A_i'(|\eta_\om^\ep(z)|) &\approx
-\frac{\big|\eta_\om^\ep(z)\big|^{1/4}}{2 \sqrt{\pi}}e^{-\frac{2}{3}
  \big|\eta_\om^\ep(z)\big|^{3/2}} \left[1 +
  O\left(\big|\eta_\om^\ep(z)\big|^{-3/2}\right)\right] ,
\end{align}
and 
\begin{align}
B_i(|\eta_\om^\ep(z)|) &\approx \frac{1}{\sqrt{\pi}
  |\eta_\om^\ep(z)|^{1/4}} e^{\frac{2}{3}
  \big|\eta_\om^\ep(z)\big|^{3/2}} \left[1 +
  O\left(\big|\eta_\om^\ep(z)\big|^{-3/2}\right)\right] ,
\\ B_i'(|\eta_\om^\ep(z)|) &\approx
\frac{\big|\eta_\om^\ep(z)\big|^{1/4}}{\sqrt{\pi}}e^{\frac{2}{3}
  \big|\eta_\om^\ep(z)\big|^{3/2}} \left[1 +
  O\left(\big|\eta_\om^\ep(z)\big|^{-3/2}\right)\right].
\end{align}
We also get from definitions \eqref{eq:res8} and \eqref{eq:res9} that 
\begin{equation}
\big|\eta_\om^\ep(z)\big|^{1/4} = \ep^{-1/6} [\mu^2(z)-k^2(\om)]^{1/4} Q_\om(z),
\end{equation}
and 
\begin{equation}
\frac{2}{3} \big|\eta_\om^\ep(z)\big|^{3/2} = \ep^{-1} |\phi_\om(z)|.
\end{equation}
The statement of the lemma follows by straightforward calculations
from these results and definitions \eqref{eq:res15}--\eqref{eq:res16}.
$\Box$
\section{Details on the diffusion limit} 
\label{ap:diflimit}
Here we derive the diffusion limit results stated in Lemmas \ref{lem.D1} and \ref{lem.D2}. 
\subsection{Single-frequency diffusion limit}
\label{ap:diflimit1}
Let us write equation \eqref{eq:Pr9} in the form 
\begin{equation}
\partial_z \psi_\om^\ep(\zeta) = \frac{1}{\ep^{1/3}}
        {\mathfrak{F}}\big(\ep^{1/3} \zeta, \nu^\ep(\zeta),
        \zeta^\ep(\zeta),\psi_\om^\ep(\zeta)\big), \label{eq:CSF1}
\end{equation}
with real-valued function ${\mathfrak{F}}$ defined pointwise by
\begin{align}
\mathfrak{F}\big(\ep^{1/3} \zeta, \nu^\ep(\zeta),
\zeta^\ep(\zeta),\psi_\om^\ep(\zeta)\big)
=& \frac{\sigma_\ep}{ \sqrt{|\zeta|}} J_{\om}^2(\ep^{1/3} \zeta)
\mu^2 \big(Z_{\om}(\ep^{1/3} \zeta)\big) \nu^\ep(\zeta) \nonumber
\\ &\times \left\{ 1 + \sin \Big[ \psi_\om^\ep(\zeta) +
  \frac{4}{3}\zeta^\ep(\zeta)\Big]\right\},
\end{align}
for the joint process $\left(\nu^\ep(z), \zeta^\ep(\zeta)\right)$
defined in \eqref{eq:Pr10b}.

The theorem in \cite[section III]{kim1996uniform} states that the
distribution of $\psi_\om^\ep(\om)$ can be described in the limit $\ep
\to 0$, in the long interval $(\zeta_{-}^\ep,\zeta_{\rm s}^\ep]$, by the
  distribution of the diffusion process with infinitesimal generator
  $\mathscr{L}_\zeta^\ep$.  This is the second-order differential operator
$\mathscr{L}_\zeta^\ep:\mathscr{C}^2 \to \mathscr{C}^0$ given by
\begin{align}
\mathscr{L}_\zeta^\ep = &
  \int_0^{\infty} \hspace{-0.1in} d \xi \, \left< \EE \left[
    \mathfrak{F}\big(\ep^{1/3} \zeta,
    \nu^\ep(\zeta),\zeta^\ep(\zeta),\psi\big)
    \right. \right. \nonumber\\& \left. \left. \times \partial_\psi \Big(
    \mathfrak{F}\big(\ep^{1/3} \zeta, \nu^\ep(\zeta + \ep^{2/3}
    \xi),\zeta^\ep(\zeta+\ep^{2/3} \xi),\psi\big)\partial_\psi \Big)
    \right]\right>_{\zeta}, \label{eq:CSF2}
\end{align}
where $\left< ~ \right>_{\zeta}$ denotes the average over the torus,
and $\mathscr{C}^q$ is the space of real-valued functions of $\psi$
with bounded and continuous derivatives up to order $q \ge 0$. 
Note that $\ep^{1/3} \zeta$ is of order one in our domain.

Recalling definition \eqref{eq:Pr10b} of $\nu^\ep$ and the
autocorrelation \eqref{eq:autocorrel} of $\nu$, we have
\begin{align*}
\EE \left[\nu^\ep(\zeta) \nu^\ep(\zeta + \ep^{2/3} \xi)\right]
&= \cR \left(\frac{Z_\om(\ep^{1/3} \zeta + \ep \xi)}{\ep} -
\frac{Z_\om(\ep^{1/3} \zeta)}{\ep} \right) \approx \cR
\left(\xi J_{\om}(\ep^{1/3}\zeta)\right),
\end{align*}
with error of order $\ep$.  The averages over the torus are
\begin{align*}
\left<\sin \left(\psi + \frac{4}{3}\zeta^\ep(\zeta)\right)\sin
\left(\psi + \frac{4}{3} \zeta^\ep(\zeta + \ep^{2/3} \xi)\right)
\right>_\zeta \approx \frac{1}{2}\cos
\left[ 2 \left(\ep^{1/3} \zeta \right)^{1/2} \xi\right],
\end{align*}
with error of order $\ep$, and similarly,
\begin{align}
\left< \sin \left(\psi + \frac{4}{3}\zeta^\ep(\zeta)\right)\cos
\left(\psi + \frac{4}{3} \zeta^\ep(\zeta + \ep^{2/3}
\xi)\right)\right>_\zeta \approx -\frac{1}{2}\sin \left[ 2
  \left(\ep^{1/3} \zeta \right)^{1/2} \xi\right],
\end{align}
and 
\begin{equation}
\left<\sin \left(\psi + \frac{4}{3}\zeta^\ep(\zeta)\right) \right>_\zeta = 0.
\end{equation}
Substituting in the expression of $\mathscr{L}_\zeta^\ep$, and
changing variables of integration we obtain
\begin{align}
\mathscr{L}_\zeta^\ep &\approx
a^\ep(\zeta) \partial^2_\psi+b^\ep(\zeta) \partial_\psi ,
\end{align}
where
\begin{align}
a^\ep(\zeta) &= 
 \frac{\sigma_\ep^2 J_\om^3(\ep^{1/3}
\zeta) \mu^4(Z_\om(\ep^{1/3}\zeta))}{2 \zeta}  \left[
    \hat \cR(0) + \frac{1}{2} \hat \cR \left(\frac{2 (\ep^{1/3}
      \zeta)^{1/2}}{J_\om(\ep^{1/3} \zeta)}\right) \right] , \\
b^\ep(\zeta) &=  -  \frac{\sigma_\ep^2 J_\om^3(\ep^{1/3}
\zeta) \mu^4(Z_\om(\ep^{1/3}\zeta))}{2 \zeta} 
\left[ \int_0^\infty d\xi \, \cR(\xi) \sin\left( \frac{2
    (\ep^{1/3} \zeta)^{1/2} \xi}{J_\om(\ep^{1/3} \zeta)}\right) \right].
\label{eq:C3}
\end{align}
Recall that $\hat \cR$ is the power spectral density, the Fourier transform of $\cR$.

Since the generator $\mathscr{L}_\zeta^\ep$ is a parabolic operator with constant
coefficients, the corresponding process is Gaussian. 
The random variable $\psi_\om(\zeta_{\rm s}^\ep)$ is asymptotically Gaussian 
with mean
$$
{\cal M}  = \lim_{\ep \to 0} \int_{\zeta_{-}^\ep}^{\zeta_{\rm s}^\ep} d \zeta \, b^\ep (\zeta),
$$
and variance
$$
{\cal V} = \lim_{\ep \to 0} 2 \int_{\zeta_{-}^\ep}^{\zeta_{\rm s}^\ep} d \zeta\, a^\ep (\zeta).
$$
Let us denote
\begin{equation}
a(\zeta) = \frac{3\sigma_\ep^2 k^4(\om)}{4 \gamma_\om \zeta}
\hat \cR(0)  ,
\end{equation}
where we recall that $J_\om(0) = \gamma_\om^{-1/3}$ and that
$\mu(Z_\om(0)) = \mu(z_T(\om)) = k(\om)$. Using the estimates
\begin{align*}
J_\om^3(\ep^{1/3} \zeta) \mu^4(Z_\om(\ep^{1/3}\zeta)) -
\frac{k^4(\om)}{\gamma_\om} &= O(\ep^{1/3} \zeta), \\
\hat \cR\left(\frac{2 (\ep^{1/3} \zeta)^{1/2}}{J_\om(\ep^{1/3}
  \zeta)}\right) - \hat \cR(0) &= O\left((\ep^{1/3}
\zeta)^{1/2}\right),
\end{align*}
and 
\begin{align*}
\sin\left( \frac{2 (\ep^{1/3} \zeta)^{1/2} \xi}{J_\om(\ep^{1/3}
  \zeta)}\right) = O\left((\ep^{1/3} \zeta)^{1/2}\right), 
\end{align*}
for $\xi$ in the support of $\cR$, we obtain using the dominated
convergence theorem and $\sigma_\ep^2 = O( |\ln \ep|^{-1})$ that 
\begin{equation}
\lim_{\ep \to 0}\left[ \int_{\zeta_{-}^\ep}^{\zeta_{\rm s}^\ep} d \zeta \, 
b^\ep(\zeta) \right] = 0, \quad \quad
\lim_{\ep \to 0}\left[ \int_{\zeta_{-}^\ep}^{\zeta_{\rm s}^\ep} d \zeta \, 
a^\ep(\zeta) -\int_{\zeta_{-}^\ep}^{\zeta_{\rm s}^\ep} d
  \zeta \, a(\zeta) \right] = 0 ,
\end{equation}
which shows that $ {\cal M} =0$ and $ {\cal V} $ is given by
\begin{equation*}
 {\cal V}  = 2 \lim_{\ep \to 0}
\int_{\zeta_{-}^\ep}^{\zeta_{\rm s}^\ep} d \zeta \, a(\zeta)
= \frac{3 k^4(\om)}{2 \gamma_\om} \hat \cR(0)
\left[ \lim_{\ep \to 0} \sigma_\ep^2 \ln
  \left(\frac{\zeta_{\rm s}^{\ep}}{\zeta_{-}^\ep}\right) \right] ,
\end{equation*}
with $\zeta_{\rm s}^\ep$ defined in \eqref{eq:zeta0} and $\zeta_{-}^\ep =
3\ep^{1/3}$. We also have  
\begin{align*}
\lim_{\ep \to 0} \sigma_\ep^2 \ln
\left(\frac{\zeta_{\rm s}^{\ep}}{\zeta_{-}^\ep}\right) &= \lim_{\ep \to 0}
\sigma_\ep^2 \ln \left[\frac{\ep^{-2/3}}{3} \left(\frac{3}{2}
  \phi_\om(0)\right)^{2/3} \right] = \frac{2}{3} \lim_{\ep \to 0}
\sigma_\ep^2 \ln \left(\frac{\phi_\om(0)}{\ep}\right),
\end{align*}
and therefore
\begin{equation}
 {\cal V}  =
\frac{k^4(\om)}{ \gamma_\om} \hat \cR(0) \left[
  \lim_{\ep \to 0} \sigma_\ep^2 \ln
  \left(\frac{\phi_\om(0)}{\ep}\right) \right]  .
\label{eq:C5}
\end{equation}
This gives the asymptotic variance of the random phase
$\psi_\om(\zeta_{\rm s}^\ep)$ at the source location $\zeta_{\rm s}^\ep$
and completes the proof of Lemma \ref{lem.D1}.
 $\Box$

\subsection{Multi-frequency diffusion limit}
\label{ap:diflimit2}
With the same argument as in section \ref{sect:diflim1}, we conclude
that we can analyze the process
\begin{equation}
\boldsymbol{\Psi^\ep}(\zeta) = \left(\psi_1^\ep(\zeta), \ldots,
\psi_m^\ep(\zeta)\right),
\end{equation}
using the diffusion limit theorem in \cite[section
  III]{kim1996uniform}. 
To apply the theorem, let us gather equations \eqref{eq:Pr16} in the
system
\begin{equation}
\partial_z \boldsymbol{\Psi}^\ep(\zeta) = \frac{1}{\ep^{1/3}}
\boldsymbol{\mathfrak{F}}\big(\ep^{1/3} \zeta, \nu^\ep(\zeta),
\zeta_1^\ep(\zeta), \ldots,
\zeta_m^\ep(\zeta),\boldsymbol{\Psi}^\ep(\zeta)\big), \label{eq:CMF1}
\end{equation}
with vector-valued function $\boldsymbol{\mathfrak{F}}$ taking values in
$\mathbb{R}^m$, with components 
\begin{align}
\mathfrak{F}_j\big(\ep^{1/3} \zeta, \nu^\ep(\zeta),
\zeta_1^\ep(\zeta), \ldots, \zeta_m^\ep(\zeta),\Psi^\ep(\zeta)\big)
= &\frac{\sigma_\ep}{ \sqrt{|\zeta|}} J_{\om_o}^2(\ep^{1/3} \zeta)
\mu^2 \big(Z_{\om_o}(\ep^{1/3} \zeta)\big) \nu^\ep(\zeta) \nonumber
\\ &\times \left\{ 1 + \sin \Big[ \psi_j^\ep(\zeta) +
  \frac{4}{3}\zeta_j^\ep(\zeta)\Big]\right\}, \label{eq:CMF2}
\end{align}
for $j = 1, \ldots, m$. Recall that $\ep^{1/3} \zeta$ is order one,
because the domain is long, of order $\ep^{-1/3}$. The system
\eqref{eq:CMF1} is driven by the joint process $\left(\nu^\ep(z),
\zeta_1^\ep(\zeta), \ldots \zeta_{m}^\ep(\zeta) \right)$ on the state
space $\mathbb{R} \times [0,3 \pi/2] \times \ldots \times [0,3
  \pi/2]$, with
\begin{equation}
\nu^\ep(z) = \nu \left(\frac{Z_{\om_o}(\ep^{1/3} \zeta)}{\ep} \right),
\end{equation}
lying in $\mathbb{R}$ and
\begin{equation}
\zeta_j^\ep(z) = \left[ \ep^{-1/3} \zeta + \ep^{-1/6} w_j
  \mathcal{K}(\ep^{1/3} \zeta)\right]^{3/2}, \qquad j = 1, \ldots, m,
\end{equation}
on the torus $[0,3 \pi/2]$. 

The infinitesimal generator $\mathcal{L}_\zeta^\ep$ is a second-order
differential operator defined on real-valued functions $f(\psi_1,
\ldots, \psi_m)$, that are twice continuously differentiable, with
bounded derivatives up to order two. It is given by 
\begin{align}
\mathscr{L}_\zeta^\ep = &\sum_{j,q=1}^m \int_0^{\infty} \hspace{-0.1in}
d \xi \, \left< \EE \left[ \mathfrak{F}_j\big(\ep^{1/3} \zeta,
  \nu^\ep(\zeta),\zeta^\ep(\zeta),\psi_1, \ldots, \psi_m\big)
  \right. \right. \nonumber\\ &\left. \left. \times \partial_{\psi_j} \Big(
  \mathfrak{F}_q\big(\ep^{1/3} \zeta, \nu^\ep(\zeta + \ep^{2/3}
  \xi),\zeta^\ep(\zeta+\ep^{2/3} \xi),\psi_1, \ldots,
  \psi_m\big)\partial_{\psi_q} \Big)
  \right]\right>_{\zeta}, \label{eq:CMF3}
\end{align}
where $\left< ~ \right>_\zeta$ denotes again the average over the torus.  
We obtain as in the previous section that
\begin{align*}
\EE \left[\nu^\ep(\zeta) \nu^\ep(\zeta + \ep^{2/3} \xi)\right]
&= \cR \left(\frac{Z_{\om_o}(\ep^{1/3} \zeta + \ep \xi)}{\ep} -
\frac{Z_{\om_o}(\ep^{1/3} \zeta)}{\ep} \right) \approx \cR
\left(\xi J_{\om_o}(\ep^{1/3}\zeta)\right),
\end{align*}
with error of order $\ep$, and 
\begin{align*}
\left<\sin \left(\psi_j + \frac{4}{3}\zeta_j^\ep(\zeta)\right)\sin
\left(\psi_q + \frac{4}{3} \zeta_q^\ep(\zeta + \ep^{2/3} \xi)\right)
\right>_\zeta \approx \frac{1}{2}\cos \left[ 2 \left(\ep^{1/3} \zeta
  \right)^{1/2} \xi\right] \delta_{jq}, \\ \left<\sin \left(\psi_j +
\frac{4}{3}\zeta_j^\ep(\zeta)\right)\cos \left(\psi_q + \frac{4}{3}
\zeta_q^\ep(\zeta + \ep^{2/3} \xi)\right) \right>_\zeta \approx
-\frac{1}{2}\sin \left[ 2 \left(\ep^{1/3} \zeta \right)^{1/2}
  \xi\right] \delta_{jq},
\end{align*}
where $\delta_{jq}$ is the Kronecker delta. Substituting in
\eqref{eq:CMF3}, and changing variables of integration, we get the
following expression of the infinitesimal generator, 
\begin{align}
\mathscr{L}_\zeta^\ep \approx &
\sum_{j,q=1}^m a_{jq}^\ep(\zeta) \partial^2_{\psi_j \psi_q} +\sum_{j=1}^m b_{j}^\ep(\zeta) \partial_{\psi_j} ,
\end{align}
where
\begin{align}
a_{jq}^\ep(\zeta) =& \frac{\sigma_\ep^2 J_{\om_o}^3(\ep^{1/3}
  \zeta) \mu^4(Z_{\om_o}(\ep^{1/3}\zeta))}{2 \zeta}
  \left[ \hat \cR(0) + \frac{1}{2} \hat \cR
  \left(\frac{2 (\ep^{1/3} \zeta)^{1/2}}{J_{\om_o}(\ep^{1/3}
    \zeta)}\right)\delta_{jq} \right] ,\\
b_j^\ep(\zeta)= &  -
\frac{\sigma_\ep^2 J_{\om_o}^3(\ep^{1/3}
  \zeta) \mu^4(Z_{\om_o}(\ep^{1/3}\zeta))}{2 \zeta}
  \int_0^\infty d\xi \, \cR(\xi)
\sin\left( \frac{2 (\ep^{1/3} \zeta)^{1/2} \xi}{J_\om(\ep^{1/3}
  \zeta)}\right)   .
\label{eq:CMF4}
\end{align}
Since the generator $\mathscr{L}_\zeta^\ep$ is a parabolic operator with constant
coefficients, the corresponding process is Gaussian. 
The random vector $\boldsymbol{\Psi^\ep}(\zeta_{\rm s}^\ep)$ is asymptotically Gaussian 
with mean $( {\cal M}_j)_{j=1}^m  $ given by 
$$
{\cal M}_j
= \lim_{\ep \to 0} \int_{\zeta_{-}^\ep}^{\zeta_{\rm s}^\ep}d \zeta \,  b_j^\ep (\zeta),
$$
and covariance matrix $({\cal V}_{jq})_{j,q=1}^m$ given by 
$$
{\cal V}_{jq} = \lim_{\ep \to 0} 2 \int_{\zeta_{-}^\ep}^{\zeta_{\rm s}^\ep} d \zeta \, a_{jq}^\ep (\zeta).
$$
Proceeding as in the previous section we find that ${\cal M}_j=0$ and 
$$
{\cal V}_{jq} =  \upsilon_{\om_o}^2  \frac{(\delta_{jq} + 2)}{3}  ,
$$
which completes the proof of Lemma \ref{lem.D2}.
 $\Box$
 %

\bibliographystyle{siam} \bibliography{TURNING}

\end{document}